\documentclass[final]{siamltex}
\usepackage{latexsym, graphicx, epsfig, amsmath, amsfonts,amssymb,bm}
\usepackage{subfigure}
\usepackage{color}

\allowdisplaybreaks

\newtheorem{remark}{Remark}[section]

\def\be{\begin{equation}}
\def\ee{\end{equation}}

\newcommand{\norm}[1]{\left\lVert#1\right\rVert}

\title{Methods to Recover Unknown Processes in Partial Differential Equations Using Data}

\author{Zhen Chen\and Kailiang Wu\and Dongbin
       Xiu\thanks{Department of Mathematics,
		The Ohio State University, Columbus, OH 43210, USA.
		{\tt Emails: chen.7168@osu.edu,} {\tt wu.3423@osu.edu,} {\tt xiu.16@osu.edu.}
	}
}

\begin{document}
\maketitle

\begin{abstract}
We study the problem of identifying unknown processes embedded in time-dependent partial differential equation (PDE) using observational data, with an application to advection-diffusion type PDE. We first conduct theoretical analysis and derive conditions to ensure the solvability of the problem. We then present a set of numerical approaches, including Galerkin type algorithm and collocation type algorithm. Analysis of the algorithms are presented, along with their implementation detail. 
The Galerkin algorithm is more suitable for practical situations, particularly those with noisy data, as it avoids using derivative/gradient data. Various numerical examples are then presented to demonstrate the performance and properties of the numerical methods.
\end{abstract}


% main text

\section{Introduction} \label{sec:intro}

Data-driven discovery or identification of unknown governing equations 
 has attracted a growing amount of attention recently, from earlier
 attempts using symbolic regression
 (\cite{bongard2007automated,schmidt2009distilling}), to more recent
 work using techniques such as Gaussian processes
 \cite{raissi2018hidden}, artificial neural networks
 \cite{raissi2017physics1,raissi2017physics2}, group sparsity
 \cite{rudy2018data}, etc.  Most of the recent work transform the
 problem into an approximation problem and develop various techniques
 to create parsimonious models \cite{tibshirani1996regression}, to discover
partial differential equations \cite{rudy2017data,
	schaeffer2017learning}, and to deal with noises in data
\cite{brunton2016discovering, schaeffer2017sparse}, corruptions in data
\cite{tran2017exact}, or limited amount of data
\cite{schaeffer2017extracting}.
Methods have also been developed in conjunction with
model selection approach
\cite{Mangan20170009}, Koopman theory \cite{brunton2017chaos}, and
Gaussian process regression \cite{raissi2017machine}, etc.
Results from approximation theory have been borrowed to justify the
use of multiple short burst of trajectories \cite{WuXiu_JCPEQ18}.
More recently, machine learning methods, particularly deep neural
networks are being investigated to aid the task of equation discovery,
for ODEs \cite{E2017,raissi2018multistep,rudy2018deep,Chen_2018,qin2018data}) and PDEs
\cite{MardtPWN_Nature18, long2017pde, KhooLuYing_2018,
  han2018solving,raissi2018deep,long2018pde,WuXiu_JCP20}.

A related class of problems is to identify unknown
parameters or processess embedded in a given system of governing
equations. This is sometimes referred to as ``system
identification''. Many efforts have been devoted to this line of
research, see, e.g.
\cite{MALENGIER2008477,Karalashvili2008,Nilssen2009,Karalashvili2011,Crews2012,SCHORSCH2013193,zhuk2016source,dam2017sparse,narasingam2018data,quade2018sparse},
and more recently, \cite{raissi2018hidden,raissi2017physics2,rudy2018data}.
The majority of the existing work focused on identification of unknown
parameters, which take constant values throughout the domain of interest.
%This setup is different from those in the classical inverse problems (cf.~\cite{belov2012inverse}) that usually consider the 
% model parameter estimation using the information on the domain
% boundary or the maps between different boundary conditions.
The focus and contribution of this paper is on the identification of
unknown processes, which are functions, embedded in a given system of governing
equations. In particular, we use advection-diffusion type of partial
differential equation (PDE) as our primary application.
%More specifically, we focus 
%The topic on system identification using observed data has been
%explored in recent articles including, which were focused on scalar constant parameters or time-dependent parameters.
%In the present paper, we are interested in identifying more general
%parameters that are unknown functions of spatial
%variables. Particular attention of our study was paid on the
%convection-diffusion type equations, for which the parameter
%identification problems were also considered in different settings
%(e.g.,
Existing work on system identification for advection-diffusion
problem include
\cite{MALENGIER2008477,Karalashvili2008,Nilssen2009,Karalashvili2011,SCHORSCH2013193,zhuk2016source}),
most of which focused on identification of constant parameters.

The technical contributions of this paper include the following. 
We first present an analysis on the uniqueness of the system identification problem for 
convection-diffusion type equations. We show that separability of the
solution is  key to guarantee uniqueness.
We then present a general numerical framework for identifying
unknown functions embedded in given governing equations using  observational data of the state variable. 
This framework is based on seeking an approximation of the unknown
functions in a properly defined finite dimensional linear space, which
can be taken conveniently as the same linear space for the
discretization of the governing equation. The identification of the
unknown processes is then conducted by minimizing 
the residual of the discretized equations in certain space-time $L^2$
norm.
Under the framework, two types of algorithms, ``collocation'' and ``Galerkin'', 
are proposed, depending on the way the residues are defined and
minimized. The Galerkin algorithm utilize weak form formulation and can avoid using information of the
derivatives/gradients of the solution states. Consequently, it is more
suitable for practical computation than the collocation algorithm,
especially when measurement data contain noises.
We remark that the proposed numerical framework and algorithms are applicable to
general classes of PDEs. Our focus on advection-diffusion type PDE in
this paper is to have a concrete model to conduct theoretical analysis.

This paper is organized as follows. After the basic problem setup in
Section \ref{sec:setup}, we present its uniqueness analysis in Section
\ref{sec:analysis}. The numerical approaches are discussed in 
Section
\ref{sec:methods}, with both a general framework and two types of
algorithms, Galerkin and collocation. 
An extensive set of numerical examples are presented in Section
\ref{sec:examples},

\section{Problem Setup} \label{sec:setup}

Let $D \subseteq  \mathbb{R}^d$, $d\geq 1$, be a spatial domain with
coordinate $\bm x = (x_1,\cdots,x_d)$, and  
$T >0$ be a real number. 
Let $u(t,{\bm x})$ be a state variable, governed by 
a time-dependent partial differential equation (PDE) 
\begin{equation}\label{geqs}
\mathcal{L}(u(t,{\bm x}); \Gamma({\bm x)}) = 0,\quad \forall (t,{\bm x}) \in [0,T] \times D,
\end{equation}
where $\mathcal L$ is an known differential operator, and $\Gamma({\bm
  x})=(
\gamma_1 (\bm x),\cdots, \gamma_P (\bm x)  )$ {\em unknown}
functions depending only on the spatial variable $\bm x$. 
Suppose observation data of the solution state $u$ are available, our goal
is to identify the unknown functions $\Gamma({\bm x})$ embedded in the governing
equation \eqref{geqs}.

In order to conduct concrete theoretical and numerical analysis, we focus on advection-diffusion type PDE
\begin{align}\label{ConDifEq}
\mathcal{L} \big(u(t,{\bm x}); \Gamma({\bm
  x})\big)
  := \frac{\partial u}{\partial t} + \nabla \cdot ({\bm \alpha({\bm
  x})} F(u)) - \nabla \cdot (\kappa({\bm x}) \nabla u)=0,
\end{align}
where   $F(\cdot)$ is flux function, $\bm \alpha(\bm x) =
(\alpha^{(1)}(\bm x),\dots,\alpha^{(d)}(\bm x) )^\top$ velocity
field and $\kappa(\bm x)$ diffusivity field. The flux $F$ is assumed
to be known, and
the unknown processes to be recovered are 
$\Gamma({\bm x}) = \{
\alpha^{(1)}(\bm x),\dots,\alpha^{(d)}(\bm x), \kappa({\bm x}) \}$.
%unknown parameter functions are
%$\Gamma = (  \alpha^{(1)}(\bm x),\dots,\alpha^{(d)}(\bm x), \kappa(\bm x) )$.
%It is worth emphasizing that the methods proposed in this paper can
%be applied to identification of parameters of general PDEs.
Throughout this paper, we assume
$\Gamma\in C^1(D)$. 
Note that even though our theoretical analysis applies to this
advection-diffusion PDE, the proposed numerical algorithms are applicable to
general type operator $\mathcal{L}$.

\section{Uniqueness Analysis}\label{sec:analysis}
%Let us consider one-dimensional spatial case and 
%Assume that $u$ is a given solution of the equation \eqref{geqs}. 
%Before addressing the numerical aspects, we would like to discuss the uniqueness problem at the PDE level of recovering the unknown parameter functions 
%from the solution $u$. 
%Specifically, in this section we analytically study that, under what conditions on the solution $u$, 
%there is a unique set $\Gamma$ of parameter functions such that the equation \eqref{geqs} is satisfied. 
%Having the uniqueness is important, otherwise the parameter identification problem would be ill-posed and not well-defined. 
%All the discussions in this section are restricted to one spatially dimensional case ($d=1$), while the 
%extension to the multidimensional cases are much more difficult and will be
%pursued in separate future work.

In this section, we present theoretical analysis
for the aforementioned recovery problem. 
%More specifically, we derive
%conditions under which the unknown process $\Gamma({\bm x}) = \{
%\alpha^{(1)}(\bm x),\dots,\alpha^{(d)}(\bm x), \kappa({\bm x}) \}$ in the
%advection-diffusion equation\eqref{ConDifEq} can be uniquely
%recovered. 
We restrict our analysis to one spatial dimension with $d=1$, as
multi-dimensional analysis becomes more challenging and
remains open. We also break down the analysis into three sub-problems:
for advection equation, for diffusion equation, and finally for advection-diffusion equation.

%In the following, we shall always assume that the unknown parameter functions to be recovered belong to $C^1(D)$. 
%Then it is natural to define the uniqueness of the unknown parameter functions as follows. 
%\begin{definition}
%	Given a solution $u(t,x)$ of the equation \eqref{geqs}, we call that 
%the unknown parameter functions to be recovered are unique if 
%% $u(x,t)$ satisfies  
%\begin{equation}\label{def}
%\# \bigg( \bigcap_{t\in [0,T]} {\mathbb A}_t \bigg) = 1,
%\end{equation}
%where ${\mathbb A}_t= \{ \Gamma: \gamma_i\in C^1(D), 1\le i \le P, {\mathcal L} ( u(t,x); \Gamma ) = 0, \forall x \in D  \}$, 
%and the notation ``$\#$'' denotes the cardinality of a set. 
%\end{definition}
%For a general PDE, it is highly challenging %to verify \eqref{def} and 
%to make an a prior judgment whether the set of unknown parameter functions is unique or not. 
%For the convection-diffusion type equations \eqref{ConDifEq}, we can derive some conditions for the desired uniqueness.    

\subsection{Advection equation}

We now consider the following linear advection equation with unknown
variable velocity field 
\begin{equation}\label{Adv1}
{\mathcal L} ( u(t,x); \Gamma(x) ) := \frac{\partial }{\partial t} u(t,x) + \frac{\partial}{\partial x} \Big( \alpha(x) F(u(t,x)) \Big) =0,
\end{equation}
where $F\in C^1(\mathbb R)$ is known and $\Gamma(x) = \alpha(x)$ is unknown.

\begin{lemma}\label{lemma:Adv1}
Let $u(t,x) \in C^1( [0,T]\times D) $ be a given solution of the equation \eqref{Adv1}. 
A sufficient and necessary condition for the uniquely determine
$\alpha({x}) \in C^1(D)$  is that: there does not exist nonzero
function $\beta(x) \in C^1(D)$ such that $\beta(x) F(u(t,x))$ is
independent of $x$. 
\end{lemma}

\begin{proof}
	We first prove the sufficiency by contradiction. Assume that 
	there is another function $\widetilde \alpha \in C^1(D)$ such that $\widetilde \alpha \not\equiv  \alpha$ and 
\begin{align*}
&
	\frac{\partial }{\partial t} u(t,x) + \frac{\partial}{\partial x} \Big( \widetilde \alpha(x) F(u(t,x)) \Big)=0.
\end{align*}
Combining it with \eqref{Adv1} 
gives 
$$
\frac{\partial}{\partial x} \Big( (\widetilde \alpha(x) - \alpha (x) ) F(u(t,x)) \Big) =0,
$$
which implies that $(\widetilde \alpha(x) - \alpha (x) ) F(u(t,x))$ does not depend on $x$. 
Therefore,  $\widetilde \alpha(x) - \alpha (x) =0,~\forall x \in D$, which leads to the contradiction. Hence the parameter function $\alpha (x)$ to be recovered is unique. 

We now prove
necessity by contradiction. Assume that there is a nonzero function $\beta(x) \in C^1(D)$ such that 
$\beta(x) F( u(t,x) )$ does not depend on $x$. Then 
$$
\frac{\partial}{\partial x}\big(\beta(x) F(u(t,x))\big) =0.
$$
This, together with \eqref{Adv1}, imply 
$$
	\frac{\partial }{\partial t} u(x,t) + \frac{\partial}{\partial x} \Big( \big(\alpha(x)+\beta(x)  \big)  F(u(t,x)) \Big)=0,
$$
which is contradictory to the uniqueness of $\alpha$.
\end{proof}

\begin{definition}
Consider a bivariate function $h:(T_1,T_2)\times \Omega \to C^1$. 
%A $C^1$-function of two variables $h(t,x)$ %\in C^1( (T_1,T_2) \times \Omega )$ 
It is called separable  if it can be written as a product of two
univariate functions
 $$
h(t,x)=f(t)g(x),\quad \forall (t,x) \in  (T_1,T_2) \times \Omega,
$$ 
where $f\in C^1(T_1,T_2)$ and $g(x) \in C^1(\Omega)$. 
%Otherwise, it will be called inseparable on $[0,T] $. 
\end{definition}

\begin{theorem}\label{thm:Adv2}
	Let $u\in C^1([0,T]\times D) $ be a given solution of the equation \eqref{Adv1}. 
	If 
        %for any given open interval $\Omega \subseteq D$, 
	%there exist a temporal interval $(T_1,T_2)\subset [0,T]$ such that 
	%$F(u(t,x))$ is not separable on $ (T_1,T_2) \times \Omega $, 
	there exists no %does not exist 
	open interval $\Omega \subseteq D$ and $(T_1,T_2)\subseteq
        (0,T]$ such that 
	$F(u(t,x))$ is separable on $(T_1,T_2) \times \Omega $,
	then $\alpha (x)$ is unique.
\end{theorem}

\begin{proof}
	We prove it by contradiction. Assume that $\alpha (x)$ is not unique. Then, according to Lemma \ref{lemma:Adv1}, there exists a nonzero function $\beta(x) \in C^1(D)$ such that $\beta(x) F(u(t,x))$ does not depend on $x$. 
	  Thus we have $\beta (x_0) \neq 0$ for some $x_0$ in the interior of $D$, and 
	  $
	  \beta(x) F(u(t,x)) = c(t),~\forall (t,x) \in [0,T]\times D
	  $, for some single-variable function $c(t) \in C^1([0,T])$. 
	  Due to the sign-preserving property of $\beta (x)$, there exists
	  an open interval $\Omega_0 \subseteq D $ containing the point $x_0$ such that 
	  $$
	  \beta (x) \neq 0, \quad \forall x \in \Omega_0.
	  $$
	  Hence 
		  $$
	F(u(t,x)) = \frac{1} {\beta(x) } c(t),~\forall (t,x) \in  [0,T]\times \Omega_0 .
	$$
	This means $F(u(t,x))$ is separable on $ [0,T] \times \Omega_0
        $, which contradicts with the assumption on $F(u(t,x))$. 
	Therefore, $\alpha (x)$ is unique.  
\end{proof}

%\begin{remark}[{Why we need inseparability?}]
%	We give a simple example to illustrate the importance of inseparability. 
%	Let us consider the recovery of parameter function 
%	$\alpha(x)$ in the convection equation \eqref{Adv1} with $F(u)=u$.
%	Assume that the given solution $u(x,t)={\rm e}^{-t}{\rm e}^{-x}$ which is separable. 
%	It is easy to verify that $\alpha(x)=-1+C {\rm e}^{x}$ is satisfied  for any constant $C$. 
%	This means there infinitely many satisfied parameter functions, and the recovery of 
%	$\alpha$ from $u$ is ill-posed in this case. 	
%\end{remark}

\subsection{Diffusion equation}

We now consider the following diffusion equation
\begin{equation}\label{Dif1}
{\mathcal L} ( u(t,x); \Gamma(x) ) := \frac{\partial }{\partial t} u(x,t)
- \frac{\partial}{\partial x} \Big( \kappa(x) \frac{\partial}{\partial
  x}u(t,x) \Big) =0, 
\end{equation}
where $\Gamma(x) = \kappa (x)$ is unknown.

\begin{lemma}\label{lemma:Dif1}
	Let $u \in C^1([0,T];C^2(D)) $ be a given solution of the equation \eqref{Dif1}. 
	A sufficient and necessary condition for the uniqueness of the
        function $\kappa({x}) \in C^1(D)$  is that: there is no
        nonzero function $\beta(x) \in C^1(D)$ such that $\beta(x)
        \frac{\partial}{\partial x} u(t,x)$ is independent of $x$. 
\end{lemma}

\begin{proof}
The proof is similar to that of Lemma \ref{lemma:Adv1} and omitted here.
\end{proof}

\begin{theorem}
	Let $u\in C^1([0,T];C^2(D)) $ be a given solution of the equation \eqref{Dif1}. 
		If for any given open interval $\Omega \subseteq D$, 
	there exist a temporal interval $(T_1,T_2)\subset [0,T]$ such that 
	$\frac{\partial} {\partial x} {u(t,x)}$ is not separable on $ (T_1,T_2) \times \Omega $, then $\kappa (x)$ is unique.
\end{theorem}

\begin{proof}
The proof is similar to that of Theorem \ref{thm:Adv2} and omitted here. 
\end{proof}

\subsection{Advection-Diffusion equation}

We now consider one-dimensional advection-diffusion equation
\begin{equation}\label{AdvDif1}
{\mathcal L} ( u(t,x); \Gamma(x) ) := \frac{\partial }{\partial t} u(t,x)
+ \frac{\partial}{\partial x} \Big( \alpha(x) u(t,x) \Big) -
\frac{\partial}{\partial x} \Big( \kappa(x) \frac{\partial}{\partial
  x}u(t,x) \Big) =0, 
\end{equation}
where $\Gamma(x) = ( \alpha(x), \kappa (x) )$ is unknown.

\begin{definition}
Consider a bivariate function $h: (T_1,T_2) \times \Omega\to C^1$.
%	A $C^1$-function of two variables $h(t,x)$ %\in C^1( (T_1,T_2)\times \Omega )$ 
It is called weakly separable if it can be written as 
$$h(t,x)=f_1(t)g_1(x)+f_2(t)g_2(x),\quad \forall (t,x) \in   (T_1,T_2)
\times \Omega,
$$ 
where $f_i\in C^1(T_1,T_2)$ and $g_i\in C^1(\Omega)$, $i=1,2$. %Otherwise, it will be called generalized inseparable on $\Omega \times (T_1,T_2) $. 
\end{definition}

%\begin{lemma}
%	If the function $H(x,t)$ is generalized inseparable on $\Omega \times [0,T] $, 
%	then both $H(x,t)$ and $\frac{\partial H(x,t)}{\partial x} $ are inseparable on $\Omega \times [0,T] $. 
%\end{lemma}

%\begin{proof}
%	Obviously, $H(x,t)$ is inseparable on $\Omega \times [0,T] $. 
%	Let us show that $\frac{\partial H(x,t)}{\partial x} $ is also inseparable by contradiction. Assume that  $\frac{\partial H(x,t)}{\partial x} $ is separable on $\Omega \times [0,T] $. Then there exist two single-variable functions 
%	$\frac{\partial H(x,t)}{\partial x} = f(t) $

%\end{proof}

\begin{theorem}
	Let $u(t,x) \in C^1([0,T];C^2(D)) $ be a given solution of the equation \eqref{AdvDif1}. 
	If 
	%for any given open interval $\Omega \subseteq D$, 
	%there exist a temporal interval $(T_1,T_2)\subset [0,T]$ such that
	there is no %does not exist 
	open interval  $\Omega \subseteq D$ and $(T_1,T_2)\subseteq
        (0,T]$ such that 
	$u(t,x)$ is weakly separable on $ (T_1,T_2) \times \Omega $, then the functions $\alpha (x)$ and $\kappa (x)$ are unique.
\end{theorem}

\begin{proof}
	 Assume that there are another two functions 
	$\widetilde \alpha (x) \in C^1(D)$ and $\widetilde \kappa (x) \in C^1(D)$ such that 
	$$
	 \frac{\partial }{\partial t} u(t,x) + \frac{\partial}{\partial x} \Big( \widetilde \alpha(x) u(t,x) \Big) - \frac{\partial}{\partial x} \Big( \widetilde \kappa(x) \frac{\partial}{\partial x}u(t,x) \Big) = 0, \quad \forall (t,x) \in   [0,T]\times D,
	$$
	which, along with \eqref{AdvDif1}, imply  
\begin{equation}\label{WKLeq00}
 \frac{\partial}{\partial x} \Big( \beta(x) u(t,x) \Big) - \frac{\partial}{\partial x} \Big( \xi(x) \frac{\partial}{\partial x}u(t,x) \Big) = 0, \quad \forall (t,x) \in  [0,T]\times D.
\end{equation}
	where $\beta:=\widetilde \alpha-\alpha \in C^1(D)$ and $\xi := \widetilde \kappa - \kappa \in C^1(D)$. 
	Note that \eqref{WKLeq00} further implies that
\begin{equation}\label{WKLeq11}
\beta(x) u(t,x)  -   \xi(x) \frac{\partial}{\partial x}u(t,x)  = c(t),  \quad \forall (t,x) \in [0,T] \times D
\end{equation}	
	for some single-variable function $c(t) \in C^1( [0,T] )$. 
		Next, we only need to show that $\beta (x)=\xi(x)=0, \forall x \in D$.

	Let us first prove $\xi(x)=0, \forall x \in D$, by contradiction. Assume that 
	$\xi \not \equiv 0$. According to the continuity of $\xi$ on $D$, 
	we have  
	$\xi(x_0)\neq 0$ for some $x_0$ belonging to the interior of $D$. 
	Due to sign-preserving property for $\xi (x)$, there exists
	an open interval $\Omega_0 \subseteq D $ containing $x_0$ such that 
	$$
	\xi (x) \neq 0, \quad \forall x \in \Omega_0.
	$$  
	Let us introduce an auxiliary positive function 
	$$
	\eta (x) := \exp \left( - \int_{x_0}^x  \frac{\beta (s)}{\xi(s)} {\rm d} s  \right) >0, \quad \forall x \in \Omega_0.
	$$
	It then follows from \eqref{WKLeq11} that 
\begin{equation*}%\label{WKLeq22}
\eta (x) \left(  -   \frac{\beta(x)}{\xi(x) } u(x,t)  +
  \frac{\partial}{\partial x}u(t,x) \right) = -  c(t) \frac{ \eta (x)
}{\xi (x)},  \quad \forall (t,x) \in  [0,T]\times \Omega_0. 
\end{equation*}		
Or, equivalently, we have  
 \begin{equation}\label{WKLeq33}
\frac{\partial}{\partial x} \Big( \eta (x)  u(t,x) \Big) = -  c(t) \frac{ \eta (x) }{\xi (x)},  \quad \forall (t,x) \in   [0,T]\times \Omega_0.
 \end{equation}	
By integrating \eqref{WKLeq33} we have 
 \begin{equation*}%\label{WKLeq44}
 \eta (x)  u(t,x) - \eta (x_0)  u(t,x_0)  = -  c(t) \int_{x_0}^x \frac{ \eta (s) }{\xi (s)} ds,  \quad \forall (t,x) \in   [0,T]\times \Omega_0.
\end{equation*}		
Note that $\eta(x_0)=1$. We then obtain
 \begin{equation}\label{WKLeq55}
  u(t,x)   =  u(t,x_0) \times \frac{1}{ \eta (x)}    +  c(t)  \times 
  \left( -\frac{1}{\eta (x) } \int_{x_0}^x \frac{ \eta (s) }{\xi (s)} ds \right),  \quad \forall (t,x) \in   [0,T] \times \Omega_0.
\end{equation}	
This implies that $u(t,x)$ is weakly separable on $ [0,T]  \times
\Omega_0 $ and contradicts with the hypothesis on $u(x,t)$. Therefore, the assumption 
	that $\xi \not \equiv 0$ is incorrect. Hence we complete the proof of $\xi(x)=0, \forall x \in D$. 
	
	By substituting $\xi(x) \equiv 0$ into \eqref{WKLeq11}, we
        then have
	\begin{equation}\label{WKLeq66}
	\beta(x) u(t,x)    = c(t),  \quad \forall (t,x) \in  [0,T] \times D.
	\end{equation}	 
	We now prove $\beta \equiv 0 $ by contradiction. Assume  $\beta \not \equiv 0 $. 
	According to the continuity of $\beta$ on $D$, 
	we have  
	$\beta(x_1)\neq 0$ for some $x_1$ belonging to the interior of $D$. 
	Due to sign-preserving property of $\beta (x)$, there exists
	an open interval $\Omega_1 \subseteq D $ containing the point $x_1$ such that 
	$$
	\beta (x) \neq 0, \quad \forall x \in \Omega_1.
	$$  
		  Hence 
	$$
	u(t,x) = \frac{1} {\beta(x) } c(t),~\forall (t,x) \in    [0,T]\times \Omega_1.
	$$
	This means $u(t,x)$ is separable, and subsequently weakly separable, on $ [0,T] \times \Omega_1 $. This is a contradiction to the hypothesis on  $u(t,x)$. 
	Therefore, the assumption that  $\beta \not \equiv 0 $ is incorrect. 
	Hence we have $\beta(x)=0, \forall x \in D$.
	
	In summary, we have proved that $\beta (x)=\xi(x)=0, \forall x \in D$. 
	In other words, $\widetilde \alpha (x) = \alpha(x)$ and  $\widetilde \kappa (x) = \kappa (x)$ for all $x\in D$. The proof is completed. 
\end{proof}

\section{Numerical Methods} \label{sec:methods}

In this section, we present our numerical methods for recovery of
unknown functions embedded in PDE by using  data of the state variables. 
We focus on the advection-diffusion type equations \eqref{ConDifEq}
discussed in the previous section, although the methods are applicable
for general PDEs.

%This section is planned in two parts. First we present two general formulations, collocation and Galerkin. Second we focus on Galerkin method and derive in details for convection-diffusion type equations.

\subsection{General Framework}

%We first present the general framework of our method, while the practical formulations and 
%implementation details will be given later. 

We seek to approximate/represent the unknown functions $\Gamma({\bm x}) = \{
\alpha^{(1)}(\bm x),\dots,\alpha^{(d)}(\bm x), \kappa({\bm x}) \}$ in 
a finite $N$-dimensional linear subspace $V_N \subset L^2(D) \cap
C^2(D)$. 
Let 
${\bm \Phi} ({\bm x}):= ( \phi_1 ({\bm x}),\dots, \phi_N ({\bm x})
)^\top$ be a basis for $V_N$.
Denote $\bm{\alpha}_N  ({\bm x}):= ({\alpha}_N^{(1)}  ({\bm x}),\dots,
{\alpha}_N ^{(d)} ({\bm x}) )^\top \in [V_N]^d $ and $\kappa_N ({\bm
  x}) \in V_N $ the finite-dimensional representation of
$\bm{\alpha}({\bm x})=
\alpha^{(1)}(\bm x),\dots,\alpha^{(d)}(\bm x))^\top$ and $\kappa({\bm
  x})$, respectively.
They can be expressed as
\begin{align}\label{eq:expansion}
{\alpha}_N^{(\ell)} ({\bm x}) = {\bm a}_\ell^\top {\bm \Phi} ({\bm x}), \ 1\le \ell \le d, \qquad 
\kappa_N ({\bm x}) = {\bm k} ^\top {\bm \Phi} ({\bm x}),
\end{align}
where the coefficient vectors $\{{\bm a}_\ell\}_{\ell=1}^d$ and ${\bm k}$ 
are to be determined.

A straightforward approach to determine these finite-dimensional
unknown functions is to minimize
the residual of $\mathcal{L}(u; \bm {\alpha}_N, \kappa_N)$
in $L^2(0,T,L^2(D))$ norm, i.e., 
\begin{equation}\label{Problem00}
\min_{  \substack{ {\bm \alpha}_N \in [V_N]^d
		\\ \kappa_N \in V_N}} \frac{1}{T} \int_0^T  \int_D \Big( \mathcal{L}(u; \bm {\alpha}_N, \kappa_N) \Big)^2  d {\bm x}  dt. 
\end{equation}
However, this minimization problem is challenging to solve, as it
involves complicated temporal and spatial integrals, as well as the
derivatives of $u$.
We now discuss how to transform this problem into a tractable one via
proper discretization.

\subsubsection{Time Discretization} 

Let $\{ t_m \}_{m=1}^M$ denote a set of time instances in $[0,T]$,
where the data of the state variable $u$ are collected. We replace the
time integral in \eqref{Problem00} by a weighted sum. Subsequently,
the optimization problem \eqref{Problem00} can be transformed into
\begin{equation}\label{Problem11}
\min_{  \substack{ {\bm \alpha}_N \in [V_N]^d
		\\ \kappa_N \in V_N}} \sum_{m=1}^M w_m   \int_D \Big( \mathcal{L} \left(u(t_m,{\bm x}); \bm {\alpha}_N({\bm x}), \kappa_N(\bm x) \right) \Big)^2  d {\bm x}, 
\end{equation}
where $\{w_m\}_{m=1}^M$ are a set of weights. Note that with a given
time instance set $\{ t_m \}_{m=1}^M$, one can choose a proper set of
weights $\{w_m\}_{m=1}^M$ such that the weighted sum in \eqref{Problem11} is a good
approximation to the time integral in \eqref{Problem00}

\subsubsection{Space Discretization} 

Upon discretization in time, we now discuss two approaches to 
simplify the spatial integral in \eqref{Problem11}. 

\begin{itemize}
\item {\textit{``Collocation'' Type Method.}} In collocation approach, we
seek to minimize \eqref{Problem11} at selected nodes in spatial
domain, i.e., at collocation points. Let $\{\bm{x}_i \}_{i=1}^{N_C}$
be such a set of nodes, we further transform  \eqref{Problem11} into
the following problem:
\begin{equation} \label{collocation}
\min_{  \substack{ {\bm \alpha}_N \in [V_N]^d
		\\ \kappa_N \in V_N}} \sum_{m=1}^M w_m \sum_{i=1}^{N_C} \bigg(  \mathcal{L}(u(t_m,{\bm x}_i); \bm {\alpha}_N({\bm x}_i), \kappa_N ({\bm x}_i) )   \bigg)^2.
\end{equation}

%\vspace{2mm}
%\noindent
\item {\textit{``Galerkin'' Type Method.}} Let $V_{N_G} \subset L^2(D) \cap C^2(D)$ be a $N_{G}$-dimensional linear subspace, and  
$\{\phi_j ({\bm x}) \}_{j=1}^{N_G}$ be an orthonormal basis of
$V_{N_G}$. We use $V_{N_G}$ as our testing space for the residual. Our
Galerkin type method then transform  \eqref{Problem11} into the
following problem:
\begin{equation}\label{Galerkin}
\min_{  \substack{ {\bm \alpha}_N \in [V_N]^d
		\\ \kappa_N \in V_N}} \sum_{m=1}^M w_m \sum_{j=1}^{N_G} \bigg( \int_D \mathcal{L}\Big(u(t_m,\bm x); \bm {\alpha}_N({\bm x}), \kappa_N(\bm x) \Big) \phi_j (\bm x)d {\bm x} \bigg)^2.
\end{equation}

\end{itemize}

\subsection{Application to Advection-Diffusion Equation \eqref{ConDifEq}}

We now discuss the detailed formulation when applying the
aforementioned approaches to the advection-diffusion equation
\eqref{ConDifEq}. The collocation approach \eqref{collocation} requires 
direct evaluations of the equation \eqref{ConDifEq} at the collocation
points. This is straightforward to implement and requires no further
discussion. On the other hand, the implementation of the Galerkin approach \eqref{Galerkin}
requires further discussion. First, we show that the Galerkin minimization
problem \eqref{Galerkin} for the advection-diffusion \eqref{ConDifEq}
can be re-written into the minimization problem for the expansion
coefficients \eqref{eq:expansion}.
%In this section, we derive a practical formulation for problem \eqref{Galerkin} in the following theorem.

%Let $\bm{a}_l,\bm{k} \in \mathbb R^{1\times N}$ denote the expansion coefficient vectors of ${\alpha}^{(l)}_N(\bm x)$ and $ {\kappa}_N(\bm{x})$, $l = 1,\dots,N$,  respectively.

\begin{theorem}\label{thm:Galerkin2}
	Let 
	%\begin{equation}\label{c}
	$
	 \bm c = 
	\begin{bmatrix}
	\bm{a}_1^\top, &
	\cdots, &
	\bm{a}_d^\top, &
	\bm{k}^\top
	\end{bmatrix}^\top.
	$
	%\end{equation}
	 Then the problem \eqref{Galerkin} for \eqref{ConDifEq} is equivalent to
	\begin{equation}\label{Galerkin2}
	\min_{  \bm{c} \in \mathbb{R}^{(d+1)N}} 
	%\int_0^T 
	\sum_{m=1}^M w_m
	\sum_{j=1}^{N_G} \Big(  {\bm E}_j(t_m)\bm{c} - b_j(t_m) \Big)^2, %w(t) dt,
	\end{equation}
	where for $1\leq j\leq N_G$ and $1\leq i \leq N$, $b_j(t) = \int_D \frac{\partial u}{\partial t} \phi_j d \bm{x}$, and
	\begin{equation}
\begin{split}
	{\bm E}_j(t) & := -\left[
	\bm A_j^{(1)}(t), \cdots, \bm A_j^{(d)}(t), -\bm K_j(t) 
	\right] \in \mathbb R^{1\times (d+1)N},\\
	\bm A_j^{(\ell)}(t) & = \Big(A_{j1}^{(\ell)}(t),\dots,A_{jN}^{(\ell)}(t)\Big), \\
	A_{ji}^{(\ell)}(t) & =  \int_{\partial D} F(u) \phi_j \phi_i n_l dS - \int_D F(u) \frac{\partial \phi_j}{\partial x_l} \phi_i d \bm{x} , \quad  1\le \ell \le d,
	\\ 
	\bm K_j(t) & = \Big(K_{j1}(t),\dots,K_{jN}(t) \Big),
	\\ 
	K_{ji}(t) & = \int_{\partial D} \phi_i \phi_j \nabla u \cdot \bm{n} dS - \int_{\partial D} u \phi_i \nabla \phi_j \cdot \bm{n} dS\\
	& \quad + \int_D u(\nabla \phi_i \cdot \nabla \phi_j + \phi_i \Delta \phi_j) d\bm{x},
\end{split}
\label{Xi}
	\end{equation}
	with $\bm n = (n_1,\dots,n_d)$ denoting the outward unit normal vector along $\partial D$.
\end{theorem}	
	
	\begin{proof}
		We only need to prove  
		\begin{equation}\label{proofpart1}
		\int_D \mathcal{L}(u; \bm {\alpha}_N, \kappa_N) \phi_j d {\bm x} 
		=   b_j(t) - {\bm E}_j(t)\bm{c}, \quad 0 \le t \le T,
		\end{equation}
		and then substitute \eqref{proofpart1} into \eqref{Galerkin} to obtain \eqref{Galerkin2}.
		
		To show \eqref{proofpart1}, we  split $\int_D
                \mathcal{L}(u; \bm {\alpha}_N, \kappa_N) \phi_j d {\bm
                  x}$ of \eqref{Galerkin} into three terms: time
                derivative term, advection term and diffusion term, as follows.
		\begin{align}\label{operatorL}
		\begin{aligned}
		\int_D \mathcal{L}(u; \bm {\alpha}_N, \kappa_N) \phi_j d {\bm x} & = 
		\int_D \left(\frac{\partial u}{\partial t} 
		+ \nabla \cdot ({\bm \alpha }_N F(u)) 
		- \nabla \cdot (\kappa_N \nabla u) \right) \phi_j d \bm x\\
		& = b_j(t)
		+ \int_D \nabla \cdot (\bm{\alpha}_N F(u)) \phi_j d\bm x \\
		& \quad - \int_D \nabla \cdot (\kappa_N  \nabla u) \phi_j d\bm x.
		\end{aligned}
		\end{align}
		By using integration-by-part for the advection and
                diffusion terms, we have
		\begin{align}\label{convection}
		\begin{aligned}
		\int_D \nabla \cdot (\bm{\alpha}_N F(u)) \phi_j d\bm{x} & = \int_{\partial D} F(u) \phi_j \bm{\alpha}_N\cdot \bm{n}dS 
		\\
		& \quad - \int_D \nabla \phi_j \cdot (\bm{\alpha}_N F(u))d \bm{x},
		\end{aligned}
		\end{align}
		and 
		\begin{align}\label{diffusion}
		\begin{aligned}
		&\int_D \nabla \cdot (\kappa_N \nabla u) \phi_j d\bm{x}
		\\ 
		& \quad = 
		\int_{\partial D} \kappa_N \phi_j \nabla u \cdot \bm{n} dS 
		- \int_D \nabla \phi_j \cdot (\kappa_N \nabla u) d\bm{x}\\ 
		& \quad = \int_{\partial D} \kappa_N \phi_j \nabla u \cdot \bm{n} dS 
		- \int_{\partial D}  u \kappa_N \nabla \phi_j \cdot \bm{n} dS + \int_D u \nabla \cdot ( \nabla \phi_j \kappa_N) d\bm{x} \\ 
		& \quad = \int_{\partial D} \kappa_N \phi_j \nabla u \cdot \bm{n} dS 
		- \int_{\partial D}  u \kappa_N \nabla \phi_j \cdot \bm{n} dS + \int_D u (\nabla \kappa_N \cdot \nabla \phi_j +  \kappa_N \Delta \phi_j) d\bm{x}.
		\end{aligned}
		\end{align}
		%We only have to evaluate the spatial derivative of $u$ on the domain boundary $\partial D$. 
		Note that in \eqref{diffusion}, 
		the spatial derivatives of $u$ are not required in the
                interior of $D$. %only needed along the domain boundary $\partial D$, avoiding 
		% the estimate or evaluation of spatial derivative data within the domain $D$. 
		Let $\bm a_\ell=: ( a^{(\ell)}_1, \dots, a^{(\ell)}_N )^\top$ and 
		$\bm k =: ( k_1, \dots, k_N )^\top$ be the
                coefficients in \eqref{eq:expansion}. We obtain
		\begin{align*}
		\begin{split}
		\int_{\partial D} F(u) \phi_j \bm{\alpha}_N \cdot \bm{n} dS & = \int_{\partial D} F(u) \phi_j \sum_{l=1}^d \alpha_N^{(l)} n_l dS \\
		& = \sum_{l=1}^{d} \sum_{i=1}^N a_i^{(l)} \int_{\partial D} F(u) \phi_j \phi_i n_l dS,
		\end{split}
		\\
		\begin{split}
		- \int_D \nabla \phi_j \cdot (\bm{\alpha_N} F(u)) d \bm{x} & = - \int_D F(u) \sum_{l=1}^d \frac{\partial \phi_j}{\partial x_l} \alpha_N^{(l)} d \bm{x} \\
		& = - \sum_{l=1}^{d} \sum_{i=1}^{N} a_i^{(l)} \int_D F(u) \frac{\partial \phi_j}{\partial x_l} \phi_i d \bm{x},
		\end{split}
		\\
		\begin{split}
		\int_{\partial D} \kappa_N \phi_j \nabla u \cdot \bm{n} dS & = \sum_{i = 1}^N k_i \int_{\partial D} \phi_i \phi_j  \nabla u \cdot \bm{n} dS,
		\end{split}
		\\
		\begin{split}
		- \int_{\partial D} u \kappa_N \nabla \phi_j \bm{n} dS & = - \sum_{i = 1}^N k_i \int_{\partial D} u \phi_i \nabla \phi_j \cdot \bm{n} dS,
		\end{split}
		\\
		\begin{split}
		\int_D u (\nabla \kappa_N \cdot \nabla \phi_j + \kappa_N \Delta \phi_j) d\bm{x} & = \sum_{i = 1}^N k_i \int_D u(   \nabla \phi_i \cdot \nabla \phi_j  + \phi_i \Delta \phi_j) d\bm{x}.
		\end{split}
		\end{align*}
		Hence, for the advection term we have
		\begin{align}
		\nonumber
		\int_D \nabla \cdot (\bm \alpha_N F(u)) \phi_j d\bm{x} & =  \sum_{l=1}^{d} \sum_{i=1}^N a_i^{(l)} \int_{\partial D} F(u) \phi_j \phi_i n_l dS - \sum_{l=1}^{d} \sum_{i=1}^{N} a_i^{(l)} \int_D F(u) \frac{\partial \phi_j}{\partial x_l} \phi_i d \bm{x}  \\ \nonumber
		& =  \sum_{l=1}^d \sum_{i=1}^N a_i^{(l)} 
		\left(\int_{\partial D} F(u) \phi_j \phi_i n_l dS - \int_D F(u) \frac{\partial \phi_j}{\partial x_l} \phi_i d \bm{x} \right)\\
		& = 
		\begin{bmatrix} 
		\bm A_j^{(1)}(t) & \cdots & \bm A_j^{(d)}(t) 
		\end{bmatrix}
		\begin{bmatrix}
		\bm{a}_1\\
		\vdots\\
		\bm{a}_d
		\end{bmatrix}, \label{convection2}
		\end{align}
		where 
		$\bm A_j^{(l)}(t) $, $1\le l \le d,~1\le j \le N_G$, are defined in \eqref{Xi}. 
		
For the diffusion term, we have
		\begin{align} \nonumber
		\int_D \nabla \cdot (\kappa_N \nabla u) \phi_j d\bm{x} & = \sum_{i = 1}^N k_i \int_{\partial D} \phi_i \phi_j \nabla u \cdot \bm{n} dS 
		- \sum_{i = 1}^N k_i \int_{\partial D} u \phi_i \nabla \phi_j \cdot \bm{n} dS \\ \nonumber
		& \quad + \sum_{i = 1}^N k_i \int_D u(\nabla \phi_i \cdot \nabla \phi_j + \phi_i \Delta \phi_j) d\bm{x}
		\\ \nonumber
		& = \sum_{i=1}^N k_i \bigg ( \int_{\partial D} \phi_i \phi_j \nabla u \cdot \bm{n} dS - \int_{\partial D} u \phi_i \nabla \phi_j \cdot \bm{n} dS
		\\ \nonumber
		& \quad + \int_D u(\nabla \phi_i \cdot \nabla \phi_j + \phi_i \Delta \phi_j) d\bm{x} \bigg )\\
		& = \bm{K}_j(t) \bm{k},
	 \label{diffusion2}
		\end{align}
		where $\bm K_j(t) $, $1\le j \le N_G$, are defined in \eqref{Xi}.
		Combining \eqref{convection2} and \eqref{diffusion2} into \eqref{operatorL} gives 
\begin{align}\nonumber
		\begin{split}
		&\int_D \nabla \cdot (\bm \alpha_N F(u)) \phi_j d\bm{x} 
		- \int_D \nabla \cdot (\kappa_N \nabla u) \phi_j d\bm{x} 
		\\
		& \quad = 		\begin{bmatrix} 
		\bm A_j^{(1)}(t) & \cdots & \bm A_j^{(d)}(t) & -\bm K_j(t)
		\end{bmatrix}
		\begin{bmatrix}
		\bm{a}_1\\
		\vdots\\
		\bm{a}_d
		\\
		\bm k
		\end{bmatrix}
		\end{split}
		\\ \label{eq:PPPw1}
		& \quad = - {\bm E}_j (t) {\bm c}
\end{align}
		for $1 \le j \le N_G$. 
		Substituting \eqref{eq:PPPw1} into \eqref{operatorL}
		gives \eqref{proofpart1}, with which \eqref{Galerkin2} follows immediately. 
		The proof is complete.
	\end{proof}

We now derive uniqueness condition for the solution 
to the minimization problem \eqref{Galerkin2}.

\begin{theorem}\label{thm:Galerkinsolution}
Let 
	\begin{equation}\label{eq:WNotations}
\bm E (t) := 
\begin{bmatrix}
\bm E_1(t)
\\
\vdots \\
\bm E_{N_G}(t)
\end{bmatrix},
\qquad
\bm{b}(t) :=
\begin{bmatrix}
b_1(t)
\\
\vdots \\
b_{N_G}(t)
\end{bmatrix},
\end{equation}
and define 
$\bm \Xi := \sum_{m=1}^M w_m \bm E (t_m)^\top \bm E (t_m)$, which is a symmetric positive semidefinite matrix.
A solution to the minimization problem \eqref{Galerkin2} satisfies
\begin{equation}\label{www111}
\bm \Xi  {\bm c} = \sum_{m=1}^M w_m \bm E (t_m)^\top \bm b (t_m). %w(t) dt.
\end{equation}
Furthermore, if the matrix $\bm \Xi$ is nonsingular, then the problem \eqref{Galerkin2} 
has a unique solution 
\begin{equation}\label{www222}
{\bm c} = \bm \Xi^{-1} \sum_{m=1}^M w_m \bm E (t_m)^\top \bm b (t_m). %w(t) dt.
\end{equation}
\end{theorem}

\begin{proof}
We immediately have
%	\begin{align*}
%	J({\bm c}) &:= \int_0^T \sum_{j=1}^N \Big(  {\bm E}_j(t)\bm{c} - b_j(t) \Big)^2 w(t) dt
%	\\
%	& = \int_0^T  \Big(  {\bm E}(t)\bm{c} - \bm b (t) \Big)^\top \Big(  {\bm E}(t)\bm{c} - \bm b (t) \Big)  w(t)  dt
%	\\
%	& = \int_0^T \Big( \bm c ^\top \bm E(t)^\top \bm E(t) \bm c - 2 \bm c^\top {\bm E}(t)^\top \bm b (t) 
%	+ \bm b(t)^\top \bm b(t) \Big) w(t) dt
%	\\
%	& = \bm c ^\top \bm \Xi \bm c - 2 \bm c^\top \int_0^T {\bm E}(t)^\top \bm b (t) w(t) dt 
%	+ \int_0^T \bm b(t)^\top \bm b(t) w(t) dt,
%	\end{align*}
	\begin{align*}
	J({\bm c}) &:= \sum_{m=1}^M w_m \sum_{j=1}^{N_G} \Big(  {\bm E}_j(t_m)\bm{c} - b_j(t_m) \Big)^2 % w(t_m) dt
	\\
	& = \sum_{m=1}^M w_m  \Big(  {\bm E}(t_m)\bm{c} - \bm b (t_m) \Big)^\top \Big(  {\bm E}(t_m)\bm{c} - \bm b (t_m) \Big)  %w(t)  dt
	\\
	& = \sum_{m=1}^M w_m \Big( \bm c ^\top \bm E(t_m)^\top \bm E(t_m) \bm c - 2 \bm c^\top {\bm E}(t_m)^\top \bm b (t_m) 
	+ \bm b(t_m)^\top \bm b(t_m) \Big) %w(t) dt
	\\
	& = \bm c ^\top \bm \Xi \bm c - 2 \bm c^\top \sum_{m=1}^M w_m {\bm E}(t_m)^\top \bm b (t_m) %w(t) dt 
	+ \sum_{m=1}^M w_m \bm b(t_m)^\top \bm b(t_m), %w(t) dt,
	\end{align*}
	which is a positive semidefinite quadratic form in the variables $\bm c$. 
	Thus, the minima of $J({\bm c})$ satisfy \eqref{www111}. 
	If the matrix $\bm \Xi$ is nonsingular, then the linear system \eqref{www111} for $\bm c$ has unique solution, 
	which is given by \eqref{www222}. 
\end{proof}

%The optimization problem \eqref{Galerkin2} and its solution \eqref{www222} 
%both involve  
%time integrals. 
%To avoid the complicated computations of continuous integrals, 
%we propose to set the weight function $w(t)$ as 
%\begin{equation}\label{w1129}
%w(t)=\frac{1}{M} \sum_{m=1}^M \delta ( t- t_m ), 
%\end{equation}
%where $\delta(t)$ is the Dirac delta function, and $\{ t_m \}_{m=1}^M$ are some i.d.d.~random points in $[0,T]$, which will be specified in the next subsection.  
%The variables $\bm E_j(t)$ and $b_j(t)$ in \eqref{Galerkin2} also  
%involve some spatial integrals on the domain $D$ or at its interface $\partial D$, as shown in the formulas \eqref{eq:AAA2} 
%and \eqref{eq:KKK2}. 
%In general, these integrals cannot be handled analytically, and numerical quadratures should be employed. 

%The optimization problem \eqref{Galerkin2} and its solution \eqref{www222} 
%both involve  
%some spatial integrals on the domain $D$ or at its interface $\partial D$, as shown in the formulas \eqref{eq:AAA2} 
%and \eqref{eq:KKK2}. 
%In general, these integrals cannot be handled analytically, and numerical quadratures should be employed. 

%It is also worth noting that the time integral does not necessarily need to be 

\begin{remark}
Note that although the collocation approach \eqref{collocation} is
straightforward to implement, we advocate the use of the Galerkin
approach \eqref{Galerkin}. This can be seen from its implementation
for the advection-diffusion equation. Using the weak form of Galerkin
and integration-by-part, the Galerkin algorithm avoids using
derivatives of the state variable in the interior of the domain. For
many pratical problems when the spatial derivatives are not directly
available and need to be estimated from data, this is preferred
because estimating derivatives can induce more numerical errors,
especially when data contain noises.
\end{remark}

%In next section, we'll present practical implementation details.
\subsection{Implementation Detail}

%In this section, we present the implementation details of 
%our method,

Assume that $\{\tau_i\}_{i=1}^{M_{\rm tot}}$ are a large set of  
time instances in $[0,T]$, where the state $u$ are measurable. %, where $M_{\rm tot} \gg 1$. 
We set $w_m=1/M$ and $t_m= \tau_{i_m}$, $1\le m \le M$, which are $M$ uniformly i.i.d.~(independent and identically distributed) random samples 
from the set $\{\tau_i\}$. 
Let $\{\bm{x}_q\}_{q=1}^Q$ be a properly selected numerical quadrature for computing the spatial integrals on $D$ and $\partial D$. 
Based on the formulations derived in Theorems \ref{thm:Galerkin2} and \ref{thm:Galerkinsolution}, the implementation of our method proceeds as follows.

%Suppose we have highly accurate data $u(t,\bm{x})$ on $[0,T] \times D$ and collecting data may introduce Gaussian noise, 

\noindent
\textbf{Step 1: Sample Data.}
Collect the data of $u(t,\bm x)$ at the points $(t_m,{\bm x}_q)$, $1\le m \le M$, $1\le q\le Q$. 
Let us denote the sampled data as $u_{m,q} := u(t_m,\bm{x}_q) + \varepsilon_{m,q}$, where $\{\varepsilon_{m,q}\}$ are possible noises. We assume that the noises are i..i.d.~random. 
% i.i.d. random observation error.

\noindent	
\textbf{Step 2: Filter Data.} 
If the data are noisy, we propose to use a filter. For each $m=1,\dots, M$, $q=1,\dots, Q$, we locally construct a polynomial function $\widetilde{f}_{m,q}(\bm x)$ in the neighborhood of $\bm{x}_q$, and obtain filtered data $\widetilde u_{m,q} := \widetilde{f}_{m,q}(\bm x_q)$. To do so, we use standard least square minimization method and sample extra data in the neighborhood of $\bm x_q$.

\noindent	
\textbf{Step 3: Estimate Derivatives.} 
We evaluate the time derivative $\frac{\partial u}{\partial t}(t_m, \bm{x}_q)$ by locally constructing polynomial function $g_{m,q}(t)$ near $t_m$. To do so, we use standard least square minimization method and sample extra data in the neighborhood of $t_m$ from $\{\tau_i \}$. We obtain time derivative estimate $\frac{\partial u}{\partial t}(t_m, \bm{x}_q) \approx g_{m,q}'(t_m)$. Similarly, for the gradients $\nabla u(t_m, \bm{x}_q),$ $\bm x_q \in \partial D$, on the domain 
boundary, we construct (local) polynomial function $f_{m,q}(\bm x)$ for each $m=1,\dots, M$, and get $\nabla u(t_m, \bm{x}_q) \approx \nabla f_{m,q}(\bm{x}_q)$.

%We now compute $\frac{\partial u}{\partial t}(t_i, \bm{x}_k)$ on $D$ and $\nabla u(t_i, \bm{x}_k)$ on $\partial D$. For time derivative, we use central difference method. For space gradient, we first construct polynomial approximation $g_i(\bm{x})$ to $u(t_i, \bm{x})$ on $\partial D$ for fixed $t_i$, and pass gradient to polynomial functions.
	
\noindent	
\textbf{Step 4: Form $\bm E(t_m)$ and $\bm b(t_m)$.} 
We compute \eqref{Xi} at $t=t_m$ by using the filtered data $\widetilde u_{i_m,q}$ on $D$ 
and the gradient estimate $\nabla u(t_m, \bm{x}_q)$ on $\partial D$
with suitable numerical quadratures. We then compute $b_j(t_m)$ using
the derivative estimate $\frac{\partial u}{\partial t}(t_m,
\bm{x}_q)$.  
The matrix $\bm E(t_m)$ and vector $\bm b(t_m)$ in
\eqref{eq:WNotations} are then formed immediately.

\noindent	
\textbf{Step 5: Compute $\bm a_\ell$ and $\bm k$.} 
Compute the symmetric positive semidefinite matrix  
$$  {\bm \Xi} = \frac{1}{M} \sum_{m=1}^M {\bm E} (t_m)^\top {\bm E} (t_m).$$
If the matrix $ {\bm \Xi}$ is nonsingular, we obtain the unique
expansion coefficient vectors $\bm a_\ell$ and $\bm k$ in
\eqref{eq:expansion} by \eqref{www222}. That is,  
$$
\begin{bmatrix}
\bm a_1
\\
\vdots \\
\bm a_d
\\
\bm k
\end{bmatrix} =
\bm c =  {\bm \Xi}^{-1} \bigg( \frac{1}{M} \sum_{m=1}^M {\bm E} (t_m)^\top {\bm b} (t_m) \bigg),
$$
which is the minimum to the least-square problem
\begin{equation*}%\label{leastsquare11}
%{\bm c}_* = \mathop {\argmin}\limits
\min_{\bm{c} \in \mathbb{R}^{(d+1)N}}
\frac{1}{M} \sum_{m=1}^M  \sum_{j=1}^{N_G} \bigg( \bm E_j(t_m)\bm{c} - b_j(t_m) \bigg)^2
=
\min_{\bm{c} \in \mathbb{R}^{(d+1)N}} \frac{1}{M} \norm{ \begin{bmatrix}
	\bm E(t_1)
	\\
	\vdots \\
	\bm E(t_M)
	\end{bmatrix}
	\bm{c} -  	\begin{bmatrix}
	\bm b (t_1)
	\\
	\vdots \\
	\bm b (t_M)
	\end{bmatrix} }^2 .
\end{equation*}

%\textbf{Step 4: Time Integration}. We now consider \eqref{Galerkin2}. In practice, we avoid direct time integration and approximate it by summation as follows.
%	\begin{align}
%	\begin{aligned}
%	\int_0^T \sum_{j=1}^N \bigg( y_j(t) - \Xi_j(t)\bm{c} \bigg)^2  dt & \approx 
%	\sum_{m=1}^M \frac{T}{M} \sum_{j=1}^N \bigg( y_j(t_m) - \Xi_j(t_m)\bm{c} \bigg)^2\\
%	& = \sum_{m=1}^M \sum_{j=1}^N \frac{T}{M} \bigg( y_j(t_m) - \Xi_j(t_m)\bm{c} \bigg)^2
%	\end{aligned}
%	\end{align}
%	Let $j$ and $m$ run through $1$ to $N$ and $1$ to $M$ respectively with one linear ordering, and introduce the following notation
%	\begin{equation}
%	\bm \Xi := 
%	\begin{bmatrix}
%	\vdots \\
%	\Xi_j(t_m) \\
%	\vdots
%	\end{bmatrix}
%	\quad
%	\bm{y} :=
%	\begin{bmatrix}
%	\vdots \\
%	y_j(t_m) \\
%	\vdots
%	\end{bmatrix}.
%	\end{equation}
%	\eqref{Galerkin2} is equivalent to the following least-square problem
%	\begin{equation}\label{leastsquare}
%	\min_{\bm{c} \in \mathbb{R}^{(d+1)N}} \norm{\bm{\Xi}\bm{c} - \bm{y}}^2.
%	\end{equation}
%
%\noindent	
%\textbf{Step 5: Solve least square problem}. It remains to evaluate entries in $\Xi$ and $\bm{y}$. We do it by quadrature rule with $u_{i,k}, \frac{\partial u}{\partial t}(t_i, \bm{x}_k),$ and $\nabla u(t_i, \bm{x}_k)$ obtained from step $1,2$ and $3$. Once evaluation is done, we solve the least square problem \eqref{leastsquare}, get coefficient $\bm{c}$ and obtain $\bm{\alpha}_N$ and $\kappa_N$.

\section{Numerical Examples}\label{sec:examples}

In this section, we present numerical examples to demonstrate the
performance of the proposed numerical methods for advection-diffusion
problem \eqref{ConDifEq}. Our examples include both 1D and 2D cases,
as well as a nonlinear Burgers' equation that does not fall into the
category of linear advection-diffusion.

%First, we give some implementation details for our examples. Let $D$ denote space domain and $[0,T]$ time interval. Here we only consider the case where $D$ is a 1-D interval or a 2-D rectangular. 

For benchmarking purpose, we use synthetic data generated by solving
known advection-diffusion equations with high resolution.
The data are then collected over a uniformly distributed time
instances in time domain and Gauss points in spatial domain, both in
the interior and along the boundary. This results in our sets of
noiseless data. To generated noisy data, we add i.i.d. Gaussian noises
$\mathcal{N}(0,\epsilon^2)$ to the clean data, where $\epsilon$ is the
noise level. 
%
%\textit{1. Data collecting.}
%We obtain clean data $u$ by first solving PDE numerically and evaluating at sample points. In all examples, we use spectral methods in space and second order ODE solvers in time, which gives $u(t, \bm{x}), t \in \{\tau_i \}_{i=1}^{M_{tot}}, \bm{x} \in D$, where $\{ \tau_i \}_{i=1}^{M_{tot}} \subset [0, T$ is the time grid in PDE solver. We let time sample $\{t_m \}_{m=1}^M$ be $M$ uniformly selected points from $\{ \tau_i \}_{i=1}^{M_{tot}}$ and space sample $\{ \bm{x}_q\}$ be Gaussian points on $D$ and $\partial D$. We obtain noisy data by adding i.i.d. Gaussian noise $\mathcal{N}(0,\epsilon^2)$ to clean data, where $\epsilon$ being noise level.
%
%To recover the velocity and diffusivity fields,
%we employ polynomial space as our approximation space and test
%space. 

We use normalized Legendre polynomials as the basis
functions. For the noisy data cases, we employ the filtering procedure
described in the previous section. In all the examples here, we built
polynomials $\widetilde{f}_{m,q}(\bm{x})$ of degree 10 using 300
noisy data drawn from the neighborhood of $\bm{x}$. 
%
%\textit{2. Data filtering.}
%If the data is corrupted, for each sample point, we sample extra local data, build local polynomial and evaluate at sample point. In all examples, we build polynomial $\widetilde{f}_{m,q}(\bm{x})$ of degree 10 using 300 samples drawn from a neighborhood of $\bm{x}$.
%
%
%\textit{3. Derivative and gradient estimation.}
These local polynomials are also used to estimate the spatial
derivatives (when required by the algorithms). The temporal
derivatives are estimated in a similar way, by first building local
polynomials  $g_{m,q}(t)$ of degree 10 using 300 neighboring data
points and then taking their derivative. For noiseless cases, all
derivatives are computed via second-order finite difference.
%
%To estimate ime derivatives and space gradients (when required), we
%use the procedure described in the previous section. For clean
%noiseless data, this can be done by straightforward second-order
%finite difference. For noisy data, we first construct local
%polynomials $f_{m,q}(\bm x)$ and $g_{m,q}(t)$, both which with degree
%10 and by using 300 local data samples, and then take derivatives of
%these polynomials.
 % by constructing local polynomial functions $g_{m,q}(t)$, $f_{m,q}(\bm x)$, and taking derivative and gradient. If the data is clean, we construct $g_{m,q}(t)$ of order 2 using data at the nearest 3 points in $\{ \tau_i\}$, which essentially gives a second order finite difference approximation to $\frac{\partial u}{\partial t}$. For space derivative, we construct local polynomial $f_{m,q}(\bm{x})$ of degree 10 using 50 uniformly sampled points for each $\bm{x}_q \in \partial D$.  If the data is corrupted, we increase polynomial order and sample number. We construct local polynomial $f_{m,q}(\bm x)$ , $g_{m,q}(t)$ both of degree 10 with 300 samples.

%\textit{4. Approximation and  testing space.}
%We use polynomial space as approximation and testing space. For basis function, we choose normalized Legendre polynomial. We denote $d-$dimensional polynomial space of degree up to $n$ by $\mathcal{P}_d^n$. 

%\textit{5. Error computation.}
%To test accuracy of recovery, we sample equal distant points in $D$ and compute relative discrete $l^2$ error. In 2-D case, we use tensor points.

The recovered velocity and diffusivity fields are evaluated over
another set of grids and then compared to the true values. We then
report the relative $\ell^2$ errors. The sets of evaluation grids are
uniform in 1D and tensor grids in 2D.

\subsection{Example 1: Advection Equation}

We first consider 1D advection equation
\begin{equation}
\frac{\partial u(t,x)}{\partial t} = -\frac{\partial}{\partial
  x}(\alpha(x)u(t,x)), \qquad (t,x)\in (0,1]\times(-4,4),
\end{equation}
%on $[0,1] \times [-1, 1]$, 
where 
\begin{equation}
\alpha(x) = \widetilde{\alpha} (1 + \delta \sin(\omega x)),
\end{equation}
with $\widetilde{\alpha }= 0.3$, $\delta = 0.2$, $\omega = \pi$.  The clean data of $u$ is obtained by solving the equation numerically with initial condition
$$u(0,x) = \frac{1}{\sqrt{2\pi \sigma^2}} {\rm e}^{-\frac{(x-\mu)^2}{2\sigma^2}},\qquad\mu = 0, \quad\sigma^2 = 0.3. $$ 
The details of the numerical solver are listed in Table \ref{table_ex1}.
\begin{table}[htbp]
	\label{table_ex1}
	\begin{center}
		\caption{PDE solver information for convection equation in Example 1.}
		\begin{tabular}{ |c|c| }
			\hline
			data $u$ time domain  & $[0,1]$  \\  
			\hline
			data $u$ space domain  & $[-4,4]$  \\ 
			%\hline
			%initial condition & $u(0,x) = \frac{1}{\sqrt{2\pi \sigma^2}} e^{-\frac{(x-\mu)^2}{2\sigma^2}},\mu = 0, \sigma^2 = 0.3 $\\
			\hline
			boundary condition & Dirichlet condition\\
			\hline 
			scheme in time  & Crank-Nicolson \\ 
			\hline
			time step size $\Delta t$ & $10^{-4}$ \\
			\hline
			scheme in space & Chebyshev collocation  \\
			\hline
			collocation number $N_{coll}$ & 100 \\
			\hline
		\end{tabular}
	\end{center}
\end{table}
Our data of $u$ are uniformly sampled 50 points in time, and over  50
Gauss points in space, along with 2 boundary points. Our goal is to
recover $\alpha(x)$. We choose polynomial space $\mathcal{P}_{1}^{50}$ as testing space.

\begin{figure}[htbp]
	\begin{center}
		\includegraphics[width=6cm]{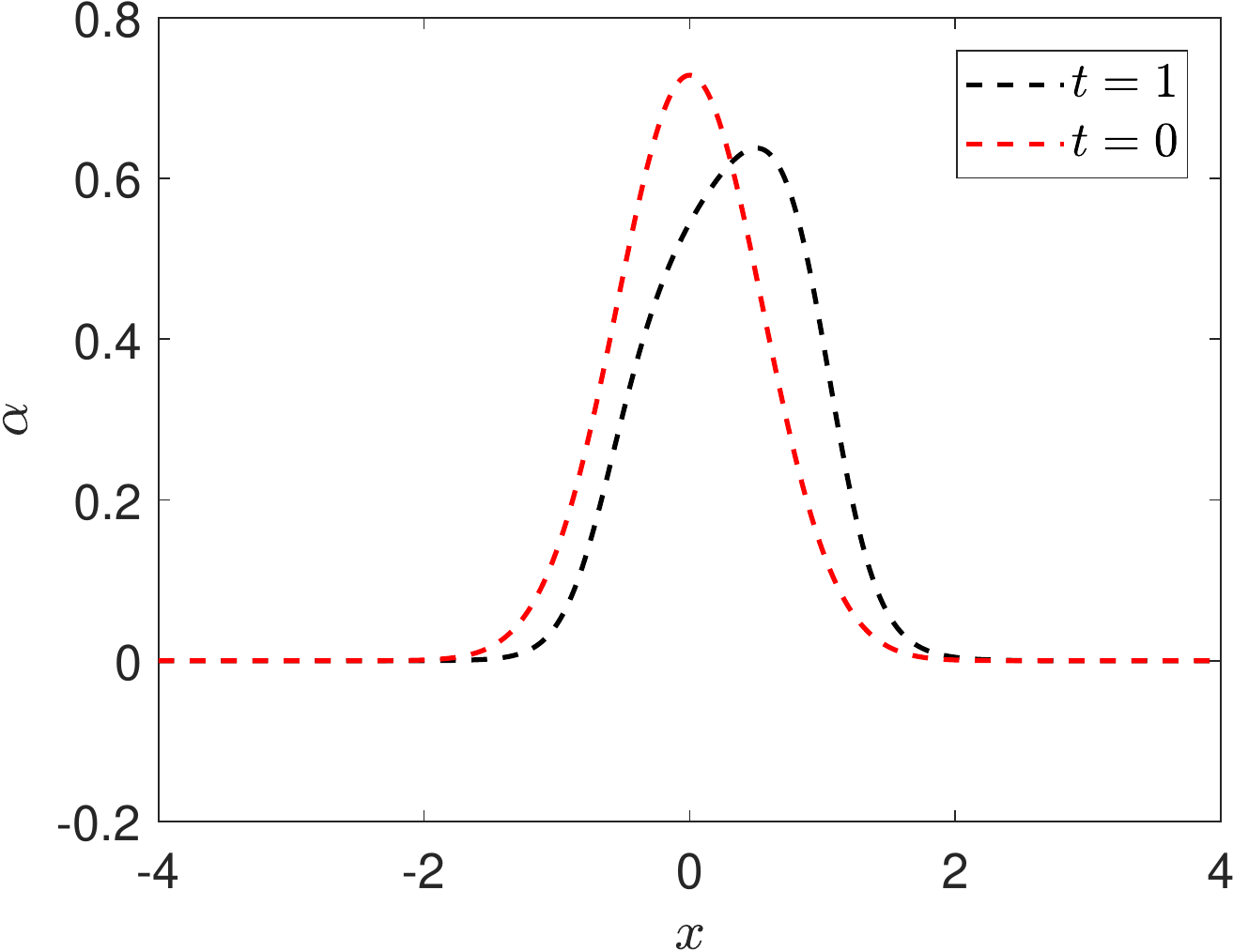}
		\includegraphics[width=6cm]{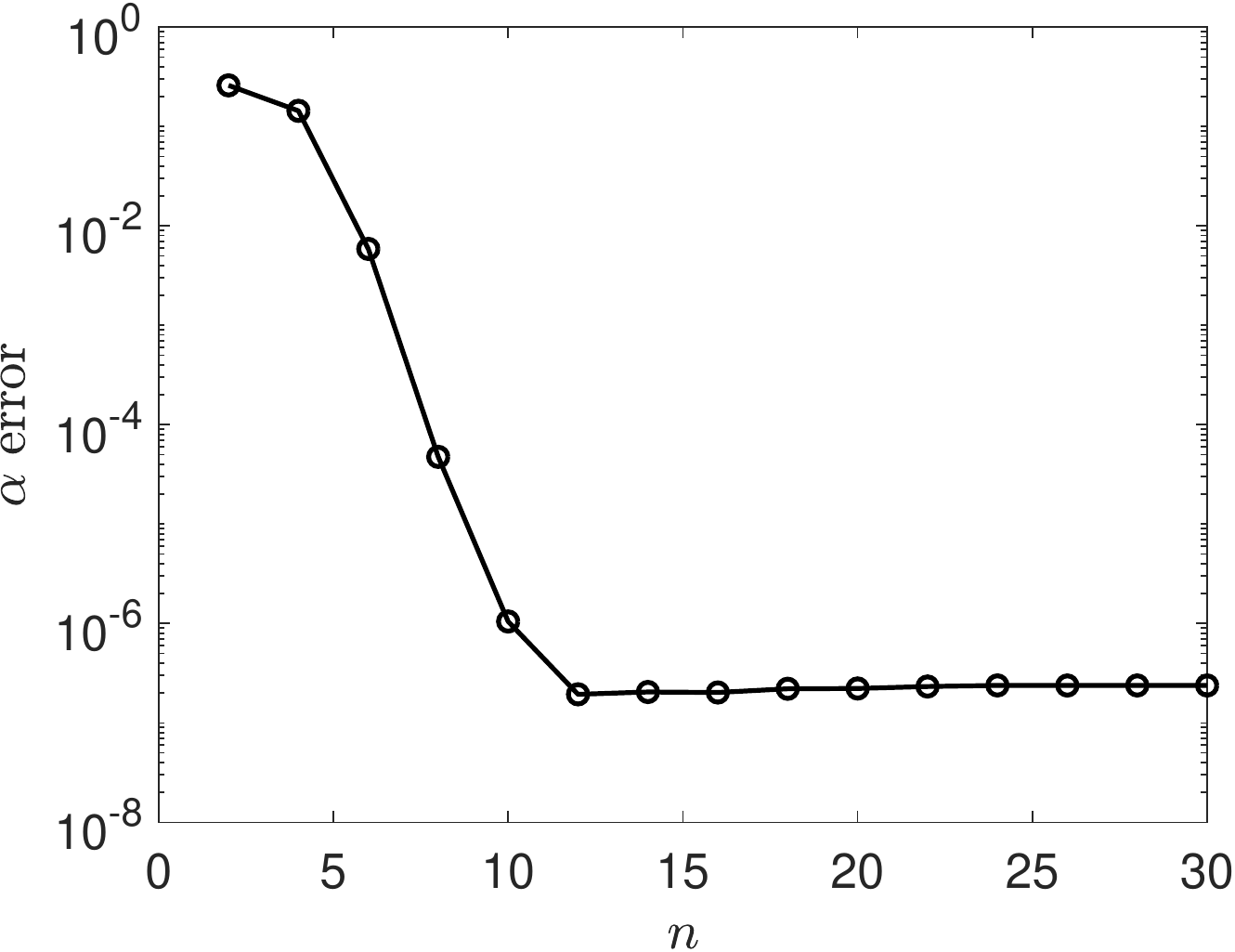}
		\caption{Example 1 with noiseless data. Left: Solution state $u$; Right:
                  Relative errors in the recovered $\alpha_n(x)$ {\em
                  vs.}  polynomial order $n$.}
		\label{fig1_clean_ex1}
\end{center}
\end{figure}

\begin{figure}
\begin{center}
		\includegraphics[width=6cm]{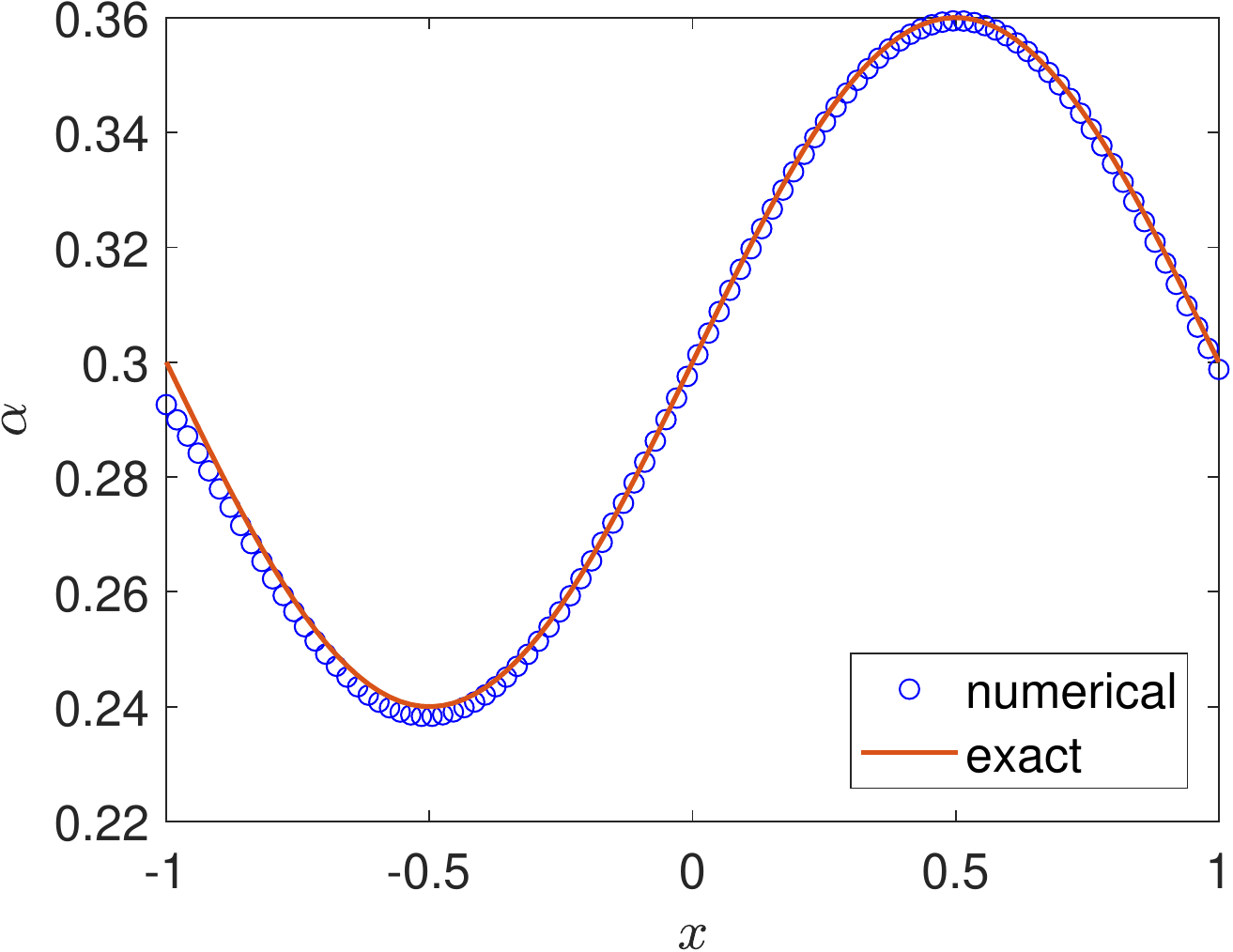}
		\includegraphics[width=6cm]{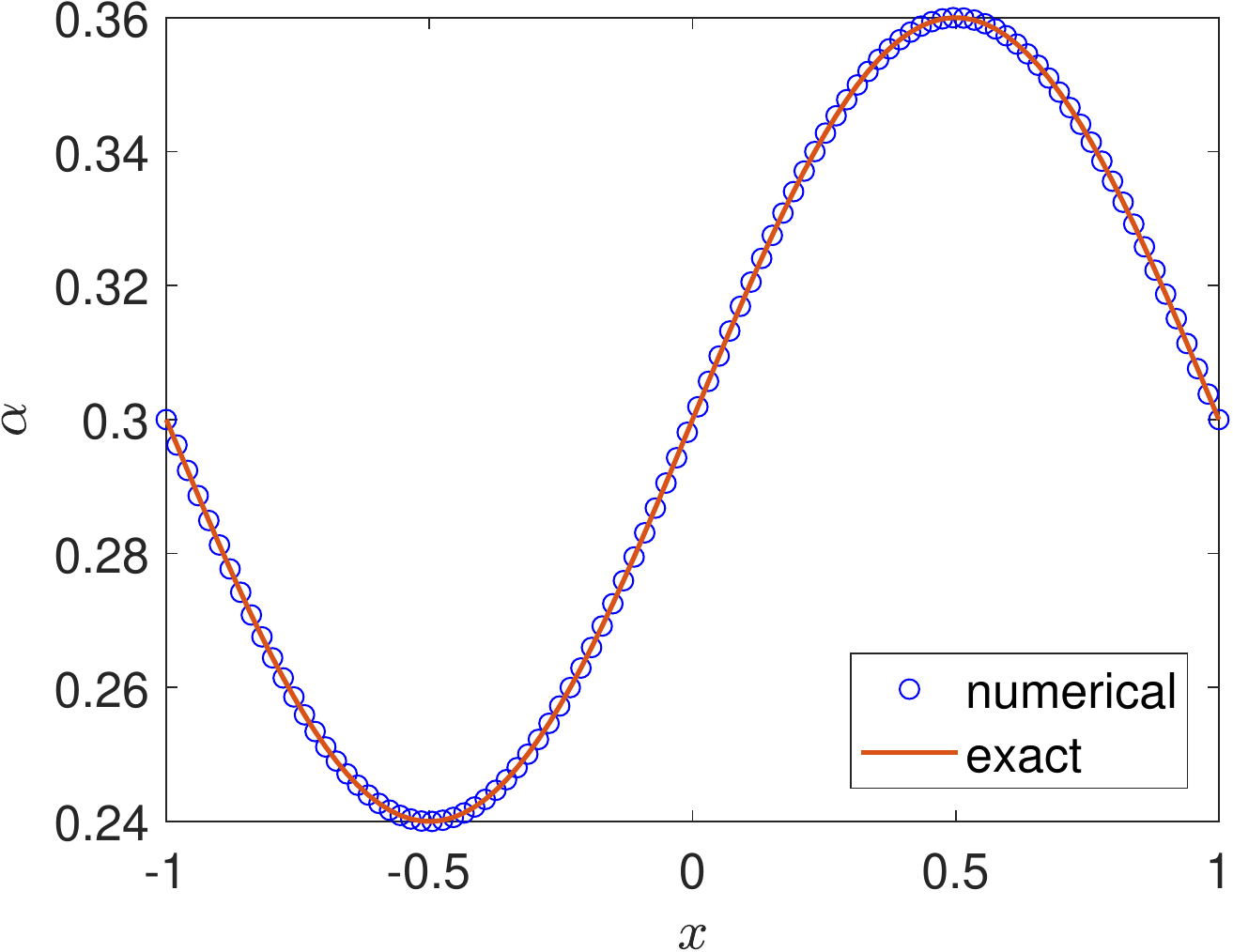}
		\caption{Example 1 with noiseless data. Left:
                  recovered $\alpha_n$  when $n=6$; Right: Recovered
                  $alpha_n(x)$ with$n=30$.}
		\label{fig1_clean_ex1_2}
	\end{center}
\end{figure}

We first consider clean noiseless data case. On the left of
Fig.~\ref{fig1_clean_ex1}, the clean data $u(0,x)$ and $u(1,x)$ are
presented. On the right of Fig.~\ref{fig1_clean_ex1}, the relative errors in our
recovered $\alpha_n(x)$ versus its polynomial order $n$ are shown.
We observe exponential decay of errors before they saturate
after $n >10$. The comparison of exact and recovered $\alpha$ are shown
in Fig.~\ref{fig1_clean_ex1_2}, for $n = 6$ and $n = 30$. 

\begin{figure}[htbp]
	\begin{center}
		\includegraphics[width=0.618\textwidth]{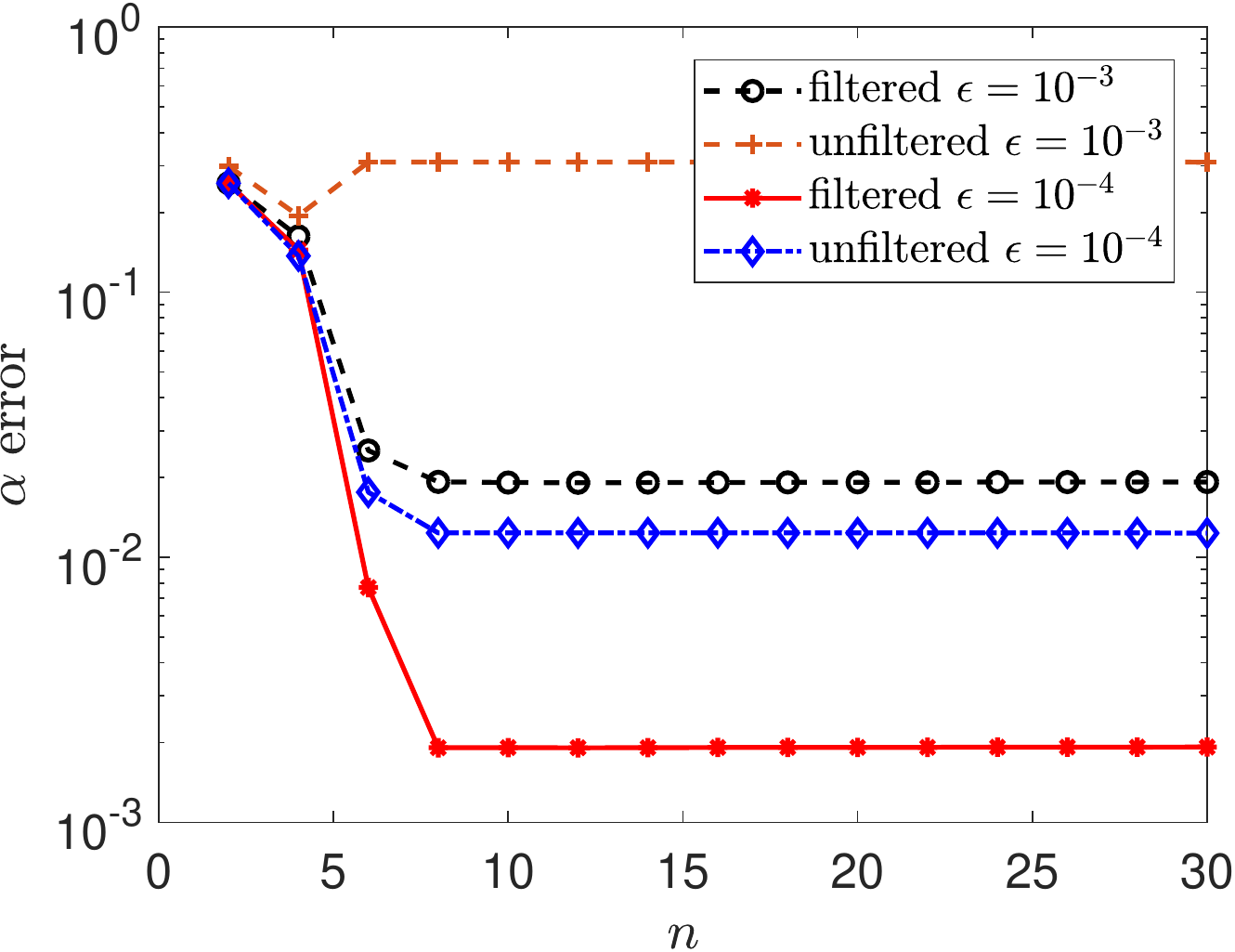}
		\caption{Example 1 with noise
                  $\mathcal{N}(0,\epsilon^2)$ in data.  Relative errors
                  in the recovered $\alpha_n(x)$ {\em vs.} polynomial
                  order $n$.}
		\label{fig2_noise_ex1}
	\end{center}
\end{figure}

Next we consider noisy data case. We add i.i.d. Gaussian noise
$\mathcal{N}(0,\epsilon^2)$ to clean data $u$, where $\epsilon =
10^{-3}$ and $10^{-4}$. The comparison of filtered and unfiltered results
is shown in Fig.~\ref{fig2_noise_ex1}. We clearly observe that
filtered results perform significantly better than the unfiltered
results, with errors one order of magnitude smaller. 
%almost $10$ times better than unfiltered one. 
This example demonstrates the necessity of employing filtering for
noisy data. 
%
%%%%%%%%%%%%%%%%%%
%
%
%
%%%%%%%%%%%%%%%%%%

\subsection{Example 2: Diffusion Equation}

We now consider a 1D diffusion equation
\begin{equation}
\frac{\partial u(t,x)}{\partial t} = \frac{\partial}{\partial x}\left(\kappa(x)\frac{\partial u(t,x)}{\partial x}\right),
\end{equation}
%on $[0,0.3] \times [-1, 1]$, 
where 
\begin{align}
\kappa(x) & = \widetilde{\kappa} \left(2 + \delta \cos(\omega x)+ 2\delta \sin \left(\frac{\omega}{2}x\right) + \delta^2 {\rm e}^x\right),
\end{align}
with $\widetilde{\kappa }= 0.3, \delta = 0.1, \omega = 4\pi$. 
The initial condition is set as
$$u(0,x) = \frac{1}{\sqrt{2\pi \sigma^2}} {\rm
  e}^{-\frac{(x-\mu)^2}{2\sigma^2}}, \quad \mu = 0, \quad \sigma^2 = 0.2.$$
The details of our numerical solver are listed in Table \ref{table_ex2}.
\begin{table}[htbp]
	\label{table_ex2}
	\begin{center}
		\caption{PDE solver information for diffusion equation in Example 2.}
		\begin{tabular}{ |c|c| }
			\hline
			data $u$ time domain  & $[0,0.3]$  \\  
			\hline
			data $u$ space domain  & $[-3,3]$   \\ 
			%\hline
			%initial condition & $u(0,x) = \frac{1}{\sqrt{2\pi \sigma^2}} e^{-\frac{(x-\mu)^2}{2\sigma^2}}, \mu = 0, \sigma^2 = 0.3 $\\
			\hline
			boundary condition & Dirichlet condition\\
			\hline 
			scheme in time  & Crank-Nicolson \\ 
			\hline
			time step size $\Delta t$ & $10^{-4}$ \\
			\hline
			scheme in space & Chebyshev collocation  \\
			\hline
			collocation number $N_{coll}$ & 150 \\
			\hline
		\end{tabular}
		
	\end{center}
\end{table}
Upon solving the equation, we collect solution data over 50 uniform points in the
temporal domain and 50 Gauss points plus 2 boundary points in the
spatial domain. We then choose $\mathcal{P}_{1}^{50}$ as the
polynomial space to recover $\kappa(x)$.

\begin{figure}[htbp]
	\begin{center}
		\includegraphics[width=6cm]{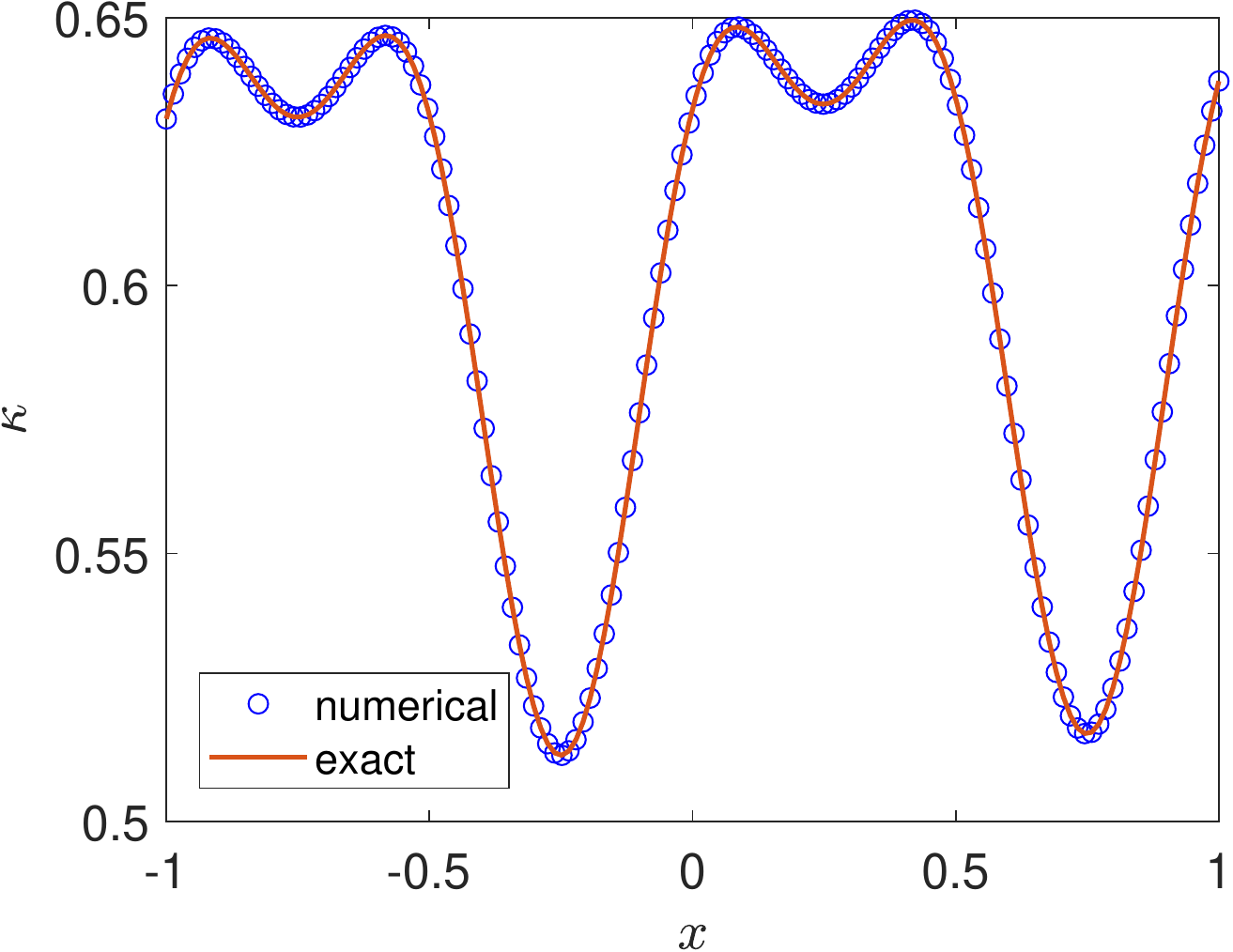}
		\includegraphics[width=6cm]{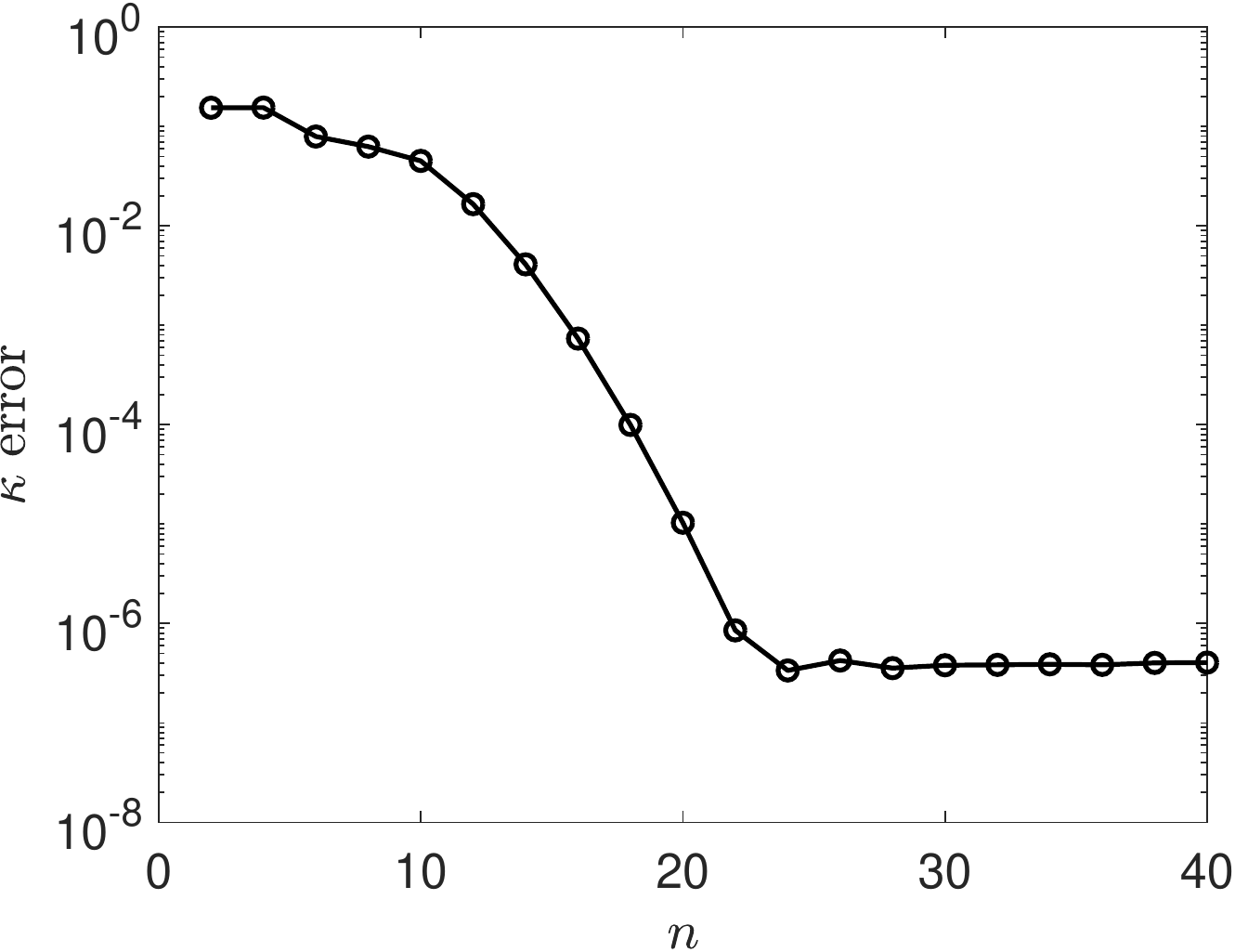}
		\caption{Example 2: Recovery of $\kappa(x)$ with
                  noiseless data. Left:
                  Result with $n = 30$; Right: Error {\em vs.}
                  polynomial order $n$.}
		\label{fig2_ex2}
\end{center}
\end{figure}

We first consider noiseless clean data case. On the left of
Fig.~\ref{fig2_ex2}, we plot the recovered $\kappa_n(x)$ with
polynomial order $n = 30$, along with the true exact $\kappa(x)$.  On the right of
Fig.\ref{fig2_ex2}, we plot error convergence and observe fast
exponential error decay. 
%Though slow at first, error decays almost exponentially when $n
%> 10$ and converges after $n = 20$. The converged error is of level
%$O(10^{-7})$. 
%

\begin{figure}[htbp]
	\begin{center}
		\includegraphics[width=0.99\textwidth]{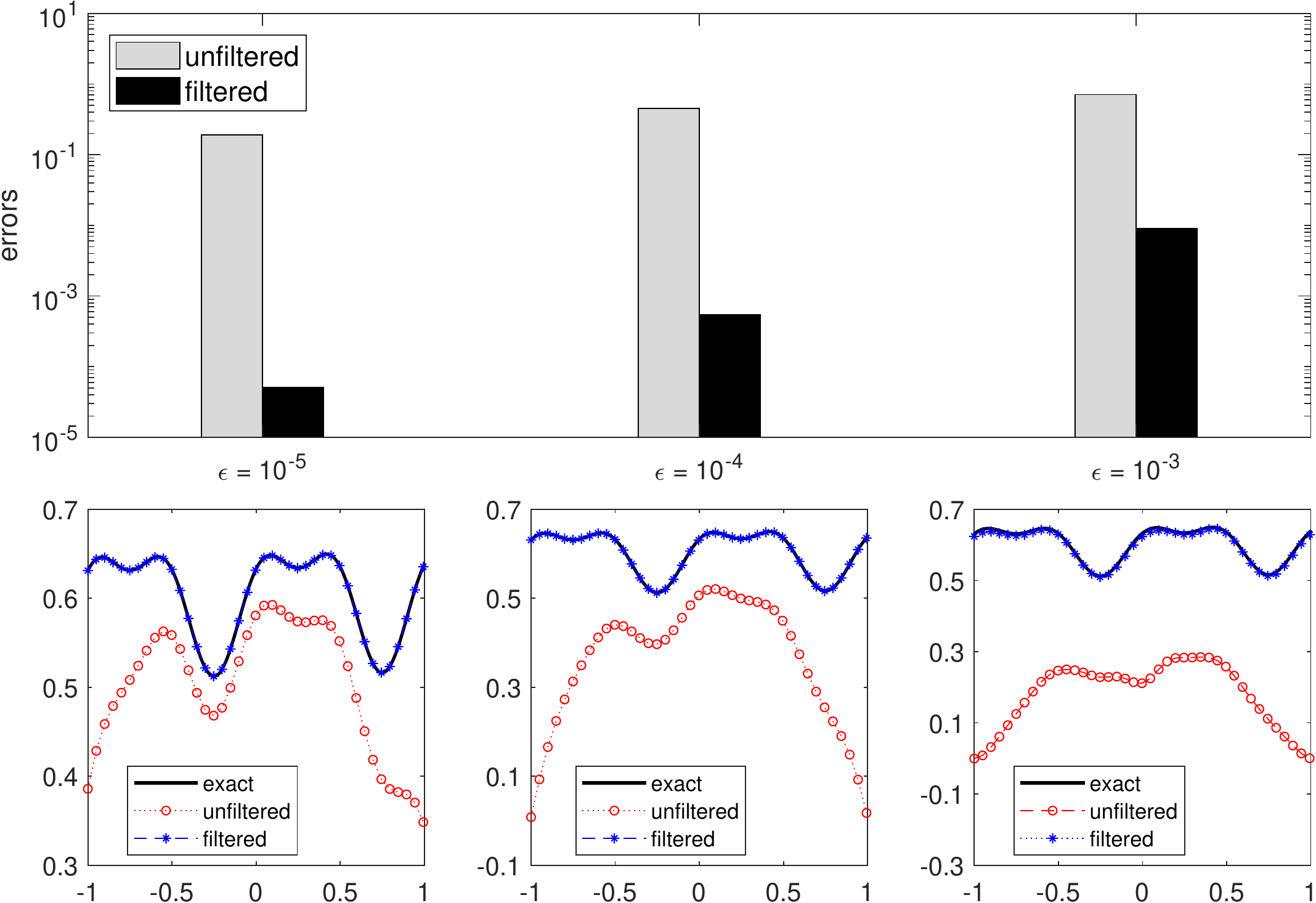}
		\caption{Example 2: Comparison of recovery results
                  with filtering and without filtering, using noisy data at different noise levels.}
		\label{fig1_ex2}
	\end{center}
\end{figure}

We then consider noisy data, with noise level at $\epsilon = 10^{-3},
10^{-4}, 10^{-5}$.
The comparison is shown
 in Fig.~\ref{fig1_ex2}, between recovery with filtering and without
 filtering. It is clearly seen that the recovery results with filtering
 are noticeably more accurate than those without filtering. 
%We can see, when filter is employed, our
 %method can capture the shape of $\kappa$ quite well. Especially in
 %the case where $\epsilon = 10^{-3}$, unfiltered solution can barely
 %recover the profile of $\kappa$ while filtered one almost get it
 %right. And though as noise level decreases, both errors decay, the
 %unfiltered error stays always above $10^{-1}$. 
%%%%%%%%%%%%%%%%%%

\begin{figure}[htbp]
\begin{center}
		\includegraphics[width=6cm]{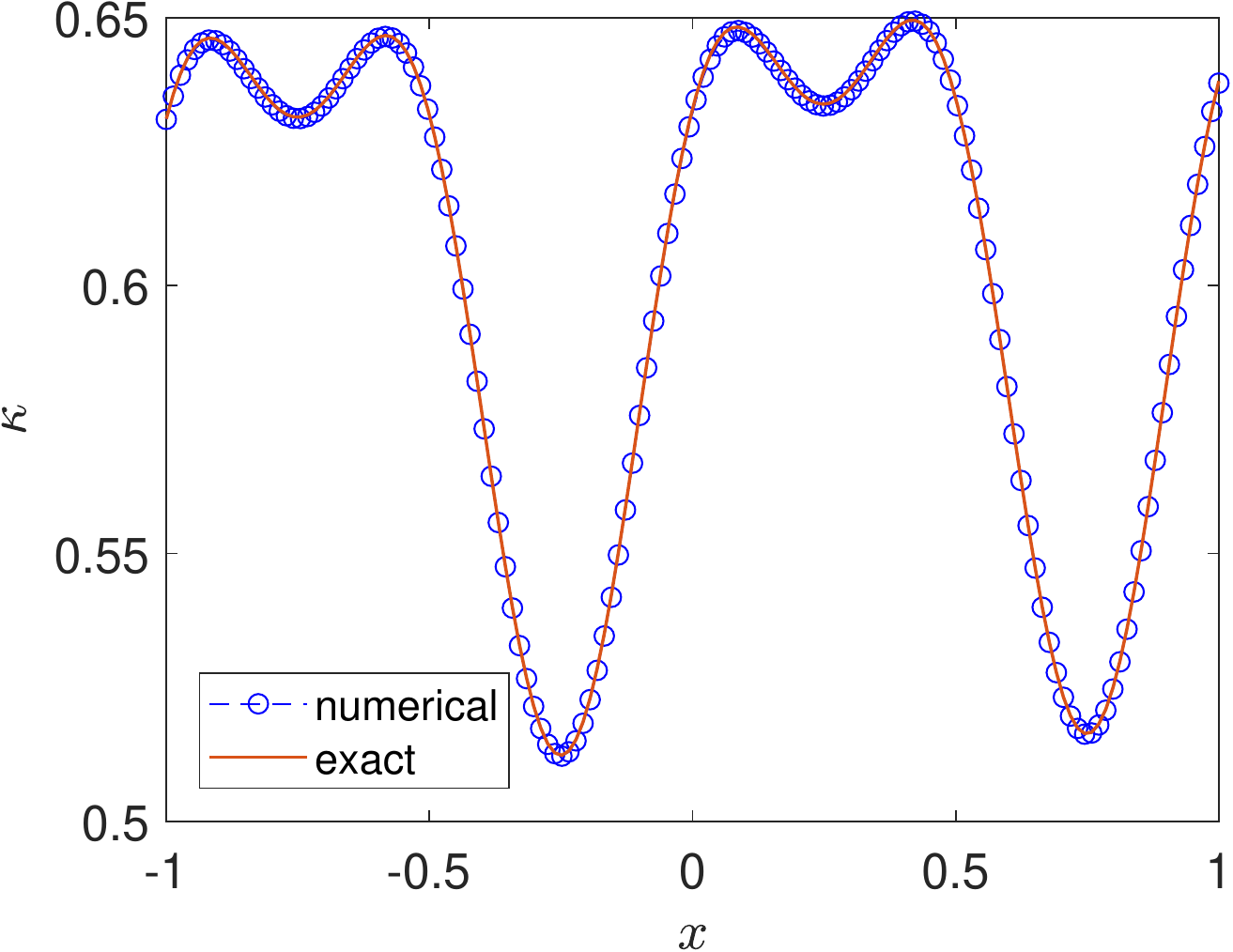}
		\includegraphics[width=6cm]{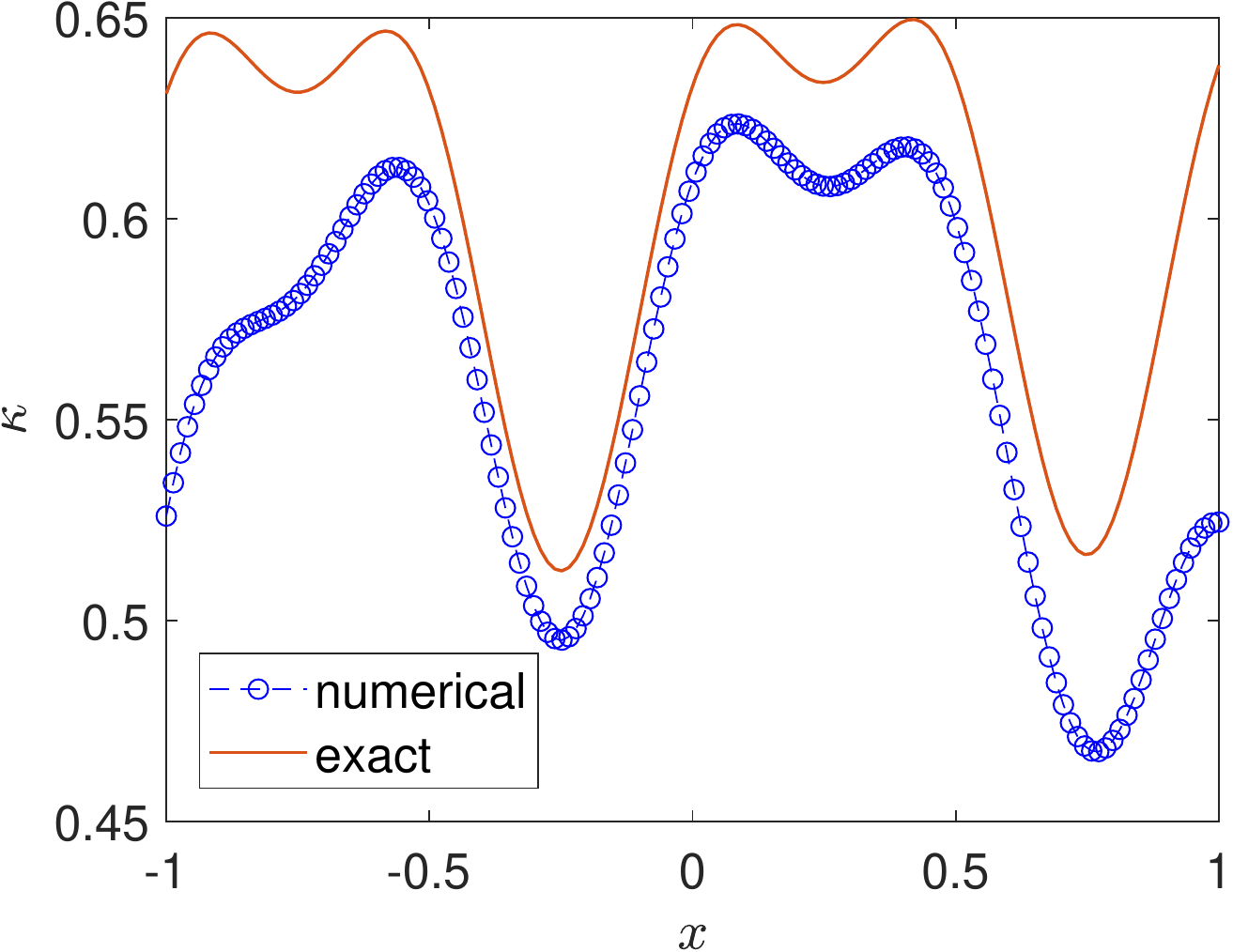}
		\caption{Example 2: $\kappa $ recovery using polynomial
                  order $n=30$ with noisy data
                  of noise level using $\epsilon = 10^{-4}$. Filtering
                  applied. Left:
                  Galerkin method; Right: Collocation method.}
		\label{fig3_ex2}	
	\end{center}
\end{figure}

\begin{figure}[htbp]
	\begin{center}
		\centering
		\includegraphics[width=0.618\textwidth]{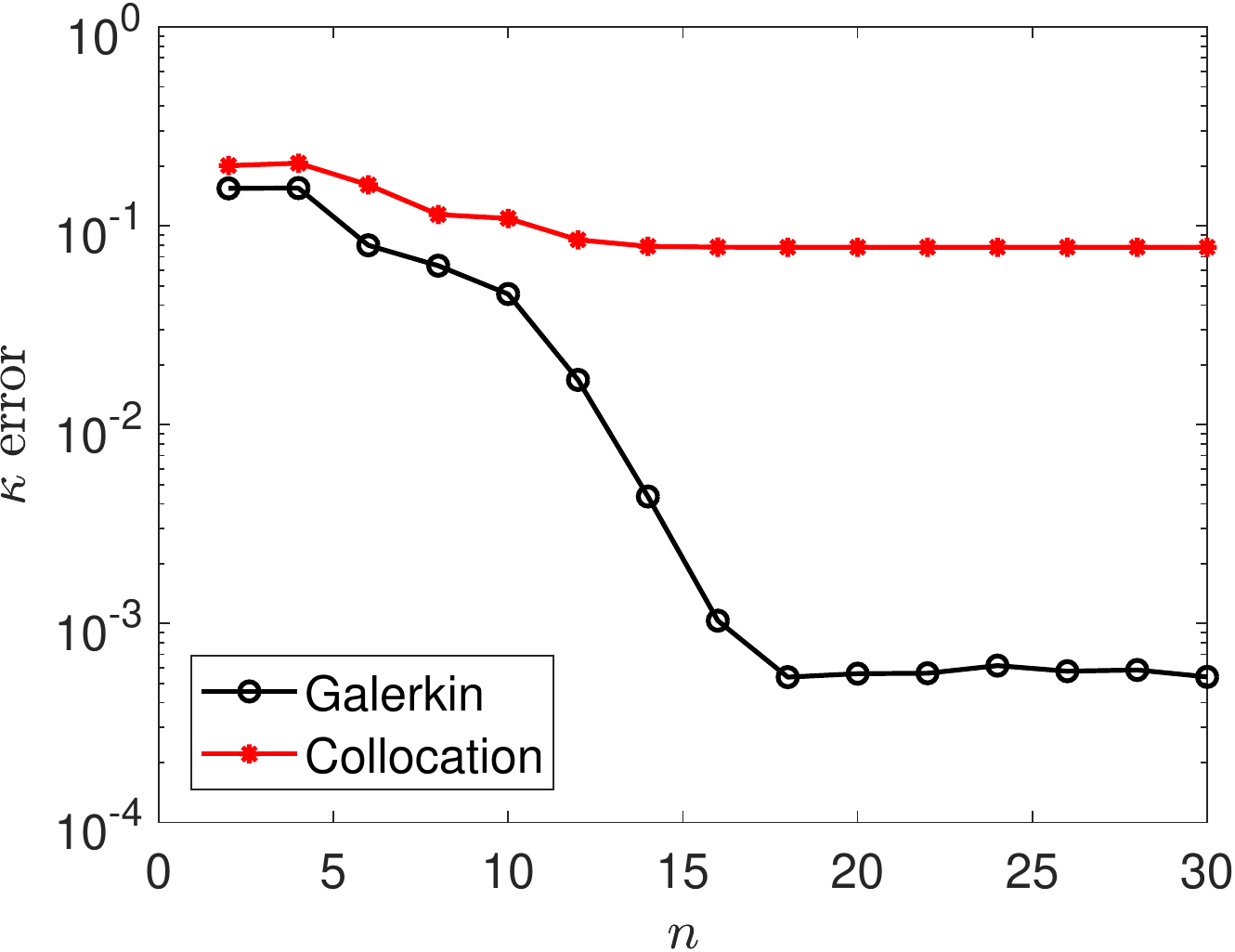}
		\caption{Example 2: noisy data with noise level
                  $\epsilon = 10^{-4}$. Filtering applied. Error {\em
                    vs.} polynomial order.}
		\label{fig4_ex2}
	\end{center}
\end{figure}

It should be mentioned that the results shown so far are obtained via
the Galerkin method. We then compare the Galerkin method and
collocation method for this example, with noisy data at noise level
$\epsilon = 10^{-4}$. Filtering is applied in both approaches. The
results are shown in Figs.~\ref{fig3_ex2} and
\ref{fig4_ex2}. Fig.~\ref{fig3_ex2} shows the results obtained with
high-order polynomial of degree $n=30$. While the Galerkin method
produces highly accurate recovery result, the results by collocation
method show visible errors and are unsatisfactory. The error convergence
with respect to increasing polynomial order is shown in
Fig.~\ref{fig4_ex2}. We can see that the collocation method fails to
converge properly as the Galerkin method does.
The primary reason for the lack of accuracy in the collocation method
is because it requires derivative estimation in the solver. Computing
derivatives with noisy data inevitably induces additional numerical
errors. On the other hand, the Galerkin method avoids much of the derivative requirement
due to its weak formulation and is able to maintain high accuracy.

\subsection{Example 3: 1D Convection-Diffusion Equation}

We now consider a 1D advection-diffusion equation
\begin{equation}
\frac{\partial u(t,x)}{\partial t} = -\frac{\partial}{\partial
  x}(\alpha(x)u(t,x)) + \frac{\partial}{\partial
  x}\left(\kappa(x)\frac{\partial u(t,x)}{\partial x}\right), \quad
\textrm{in } (0,1] \times (-1, 1),
\end{equation}
where 
\begin{equation*}
\begin{split}
\alpha(x) & = \widetilde{\alpha} \left (1 + \delta \sin(\omega x) + 2\delta \cos(\frac{\omega}{2}x) \right),\\
\kappa(x) & = \widetilde{\kappa} \left(1 + \delta \cos(\omega x)
\right),
\end{split}
\end{equation*}
with $\widetilde{\alpha} = 1$, $\widetilde{\kappa }= 0.5$, $\delta = 0.2$, $\omega = 10\pi$. 
The data set is obtained by solving the equation numerically 
with initial condition 
$$u(0,x) = \sin(\pi x) - 2{\rm e}^{-100(x-0.5)^2 } + {\rm
  e}^{-100(x+0.5)^2}. $$
Details of the numerical solver are listed in Table
\ref{table_ex3}. Data are collected over 50 uniform grids in the
temporal domain and 200 Gauss point plus the 2 boundary points in the
spatial domain. Our goal is to recover the velocity field $\alpha(x)$
and the diffusivity field $\kappa(x)$. We use $\mathcal{P}_{1}^{60}$ as
approximation and testing space. 
\begin{table}[htbp]
	\label{table_ex3}
	\begin{center}
		\caption{PDE solver information for convection-diffusion equation in Example 3.}
		\begin{tabular}{ |c|c| }
			\hline
			data $u$ time domain  & $[0,1]$  \\  
			\hline
			data $u$ space domain  & $[-1,1]$   \\ 
			%\hline
			%initial condition & \\
			\hline
			boundary condition & periodic condition\\
			\hline 
			scheme in time  & Crank-Nicolson \\ 
			\hline
			time step size $\Delta t$ & $10^{-5}$ \\
			\hline
			scheme in space & Fourier collocation  \\
			\hline
			collocation number $N_{coll}$ & 200 \\
			\hline
		\end{tabular}
	\end{center}
\end{table}

\begin{figure}[htbp]
	\begin{center}
		\includegraphics[width=0.48\textwidth]{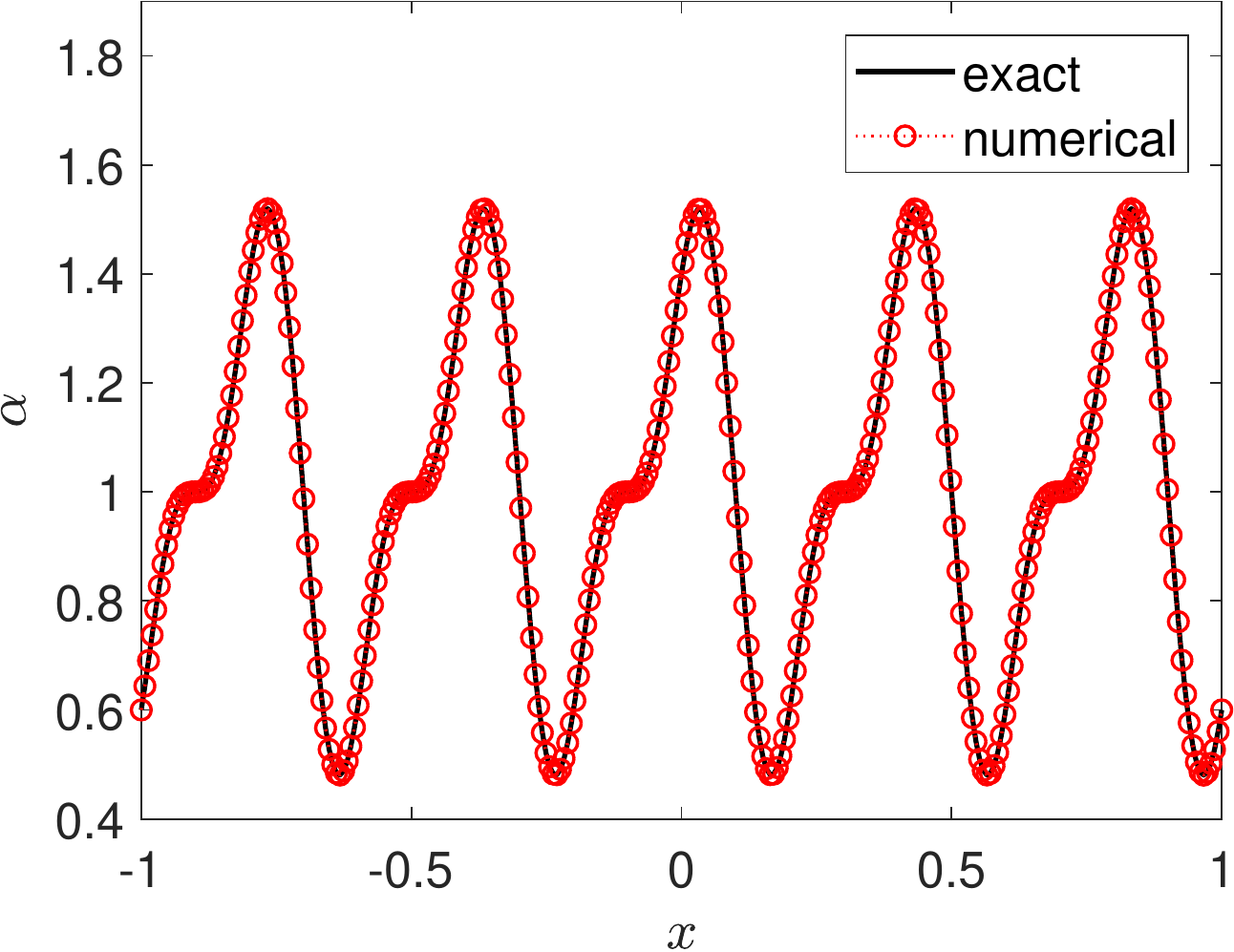}
		\includegraphics[width=0.49\textwidth]{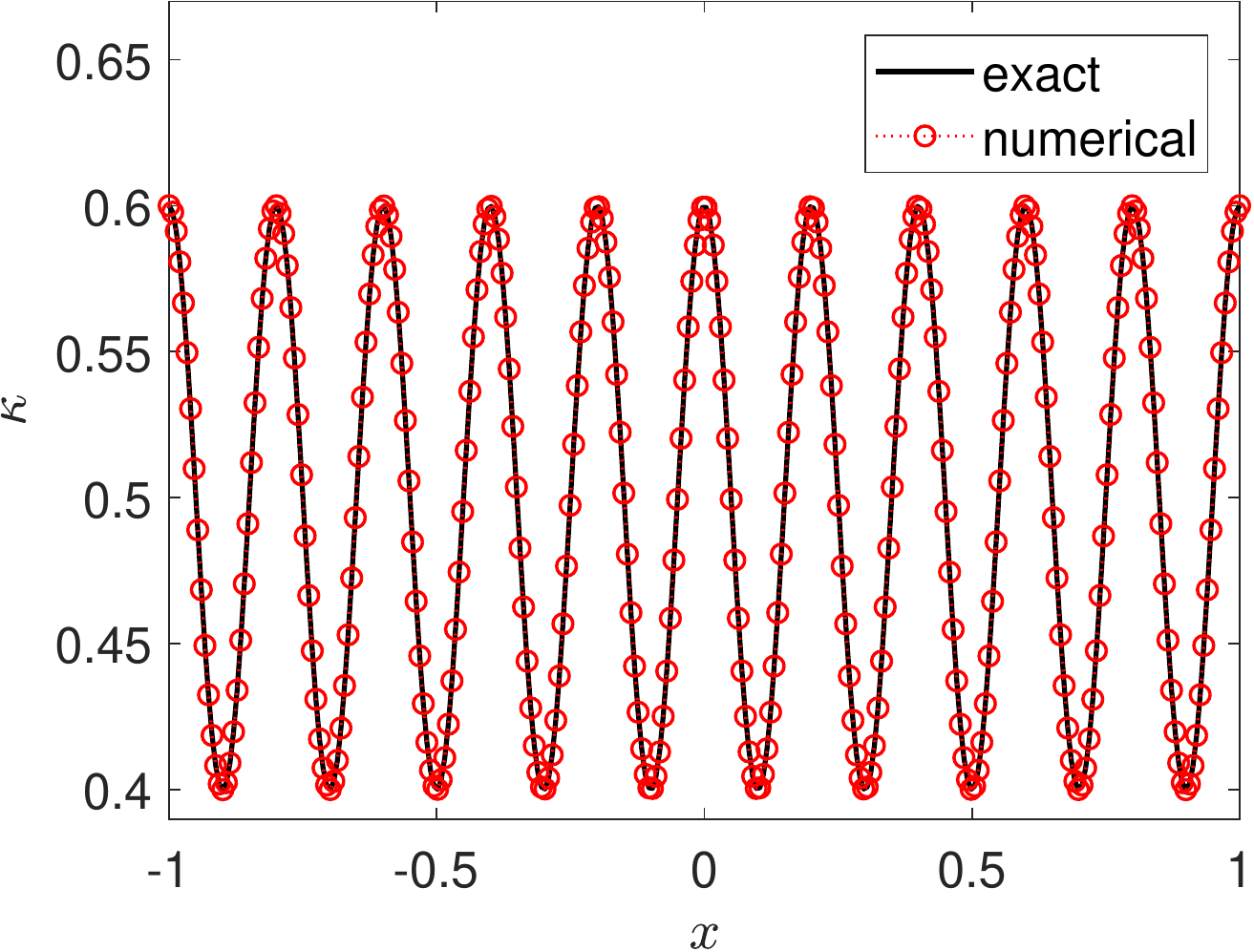}
		\caption{Example 3 with noiseless data. Left:
                  recovered $\alpha(x)$; Right: recovered $\kappa(x)$.}
		\label{fig2_ex3}
	\end{center}
\end{figure}

The recovery result for noiseless data is shown in
Fig.~\ref{fig2_ex3}. We observer excellent visual agreement between
the recovered $\alpha(x), \kappa(x)$ and their true
counterparts. Closer examination reveals that the relative errors are
$2.9335 \times 10^{-8}$ for $\alpha(x)$ and  $3.4908\times 10^{-8}$
for $\kappa(x)$. The results obtained with noisy data are not shown,
as they are visually similar to the noiseless case
and with errors dominated by the input data noise.
%%%%%%%%%%%%%%%%%%
%\begin{figure}[htbp]
%	\centering
%	\includegraphics[width=0.618\textwidth]{figures/ex3_convdiff/ex3_convdiff_solution_u_new-eps-converted-to.pdf}
%	\caption{Clean data $u$ at $t=0$ and $t=1$, for example 3.}
%	\label{fig1_ex3}
%\end{figure}

\subsection{Example 4: Viscous Burgers' Equation}

We now consider the 1D viscous Burgers' equation
\begin{equation}
\frac{\partial u(t,x)}{\partial t} = -\frac{\partial}{\partial x}\left(
\alpha(x)\frac{u(t,x)^2}{2} \right) + \frac{\partial}{\partial x}\left(
\kappa(x)\frac{\partial u(t,x)}{\partial x} \right), \quad \textrm{in } (0,0.2] \times (-1, 1),
\end{equation}
where 
\begin{equation*}
\alpha(x) = 1, \qquad
\kappa(x) = \widetilde{\kappa} \left(1 + \delta \cos(\omega x)
\right),
\end{equation*}
with $\widetilde{\kappa }= 0.1$, $\delta = 0.2$, and $\omega = 3\pi$.
The initial condition is set as
$u(0,x) = -\sin(\pi x).$

%\begin{table}[htbp]
%	\label{table_ex4}
%	\begin{center}
%		\caption{PDE solver information for Burgers' equation in Example 4.}
%		\begin{tabular}{ |c|c| }
%			\hline
%			data $u$ time domain  & $[0,0.2]$  \\  
%			\hline
%			data $u$ space domain  & $[-1,1]$   \\ 
%			%\hline
%			%initial condition & $u(0,x) = sin(x)$\\
%			\hline
%			boundary condition & periodic condition\\
%			\hline 
%			scheme in time  & Adam-Bashforth \\ 
%			\hline
%			time step size $\Delta t$ & $10^{-3}$ \\
%			\hline
%			scheme in space & Fourier collocation  \\
%			\hline
%			collocation number $N_{coll}$ & 100 \\
%			\hline
%		\end{tabular}
%	\end{center}
%\end{table}

\begin{figure}[htbp]
	\begin{center}
		\includegraphics[width=0.49\textwidth]{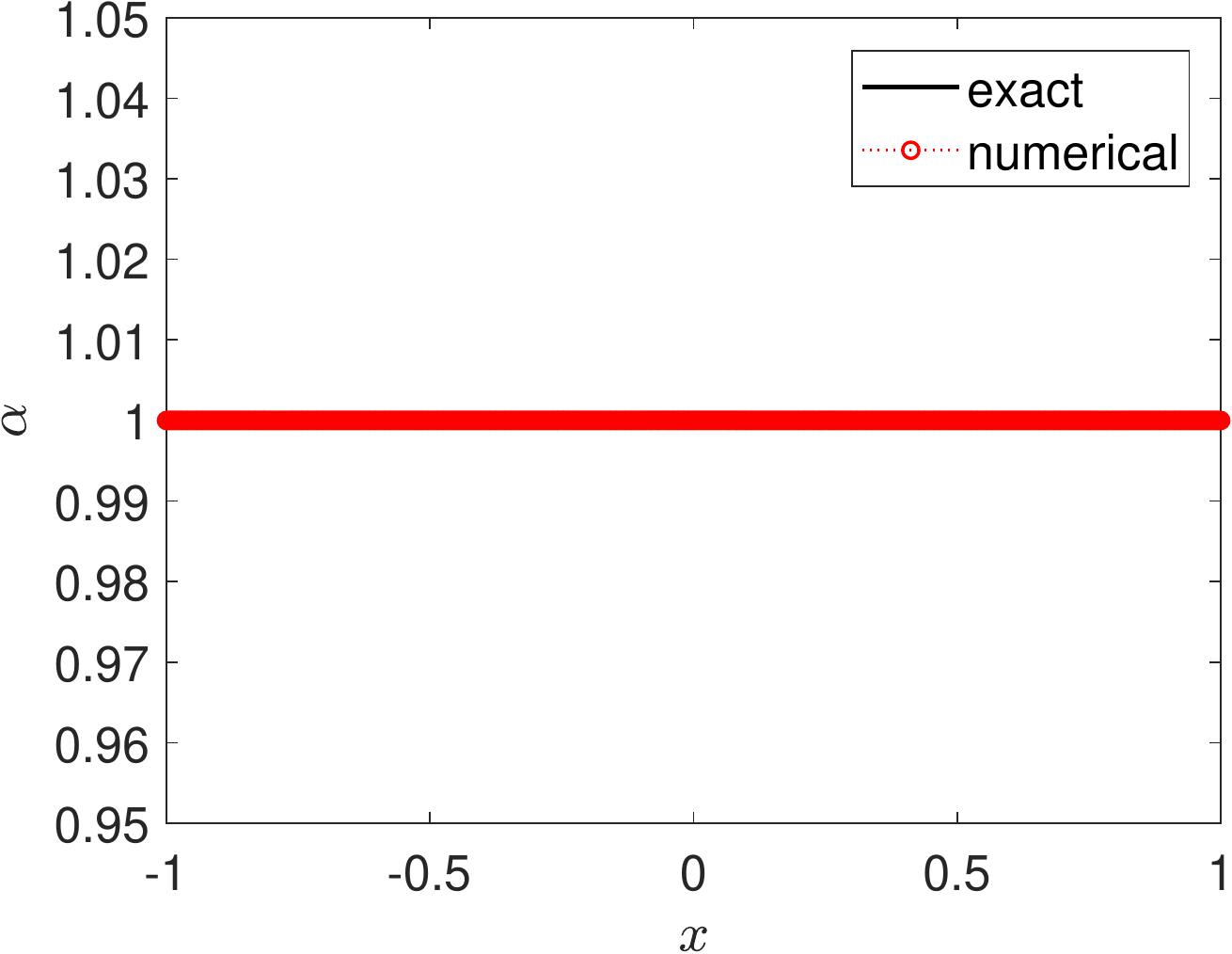}
		\includegraphics[width=0.49\textwidth]{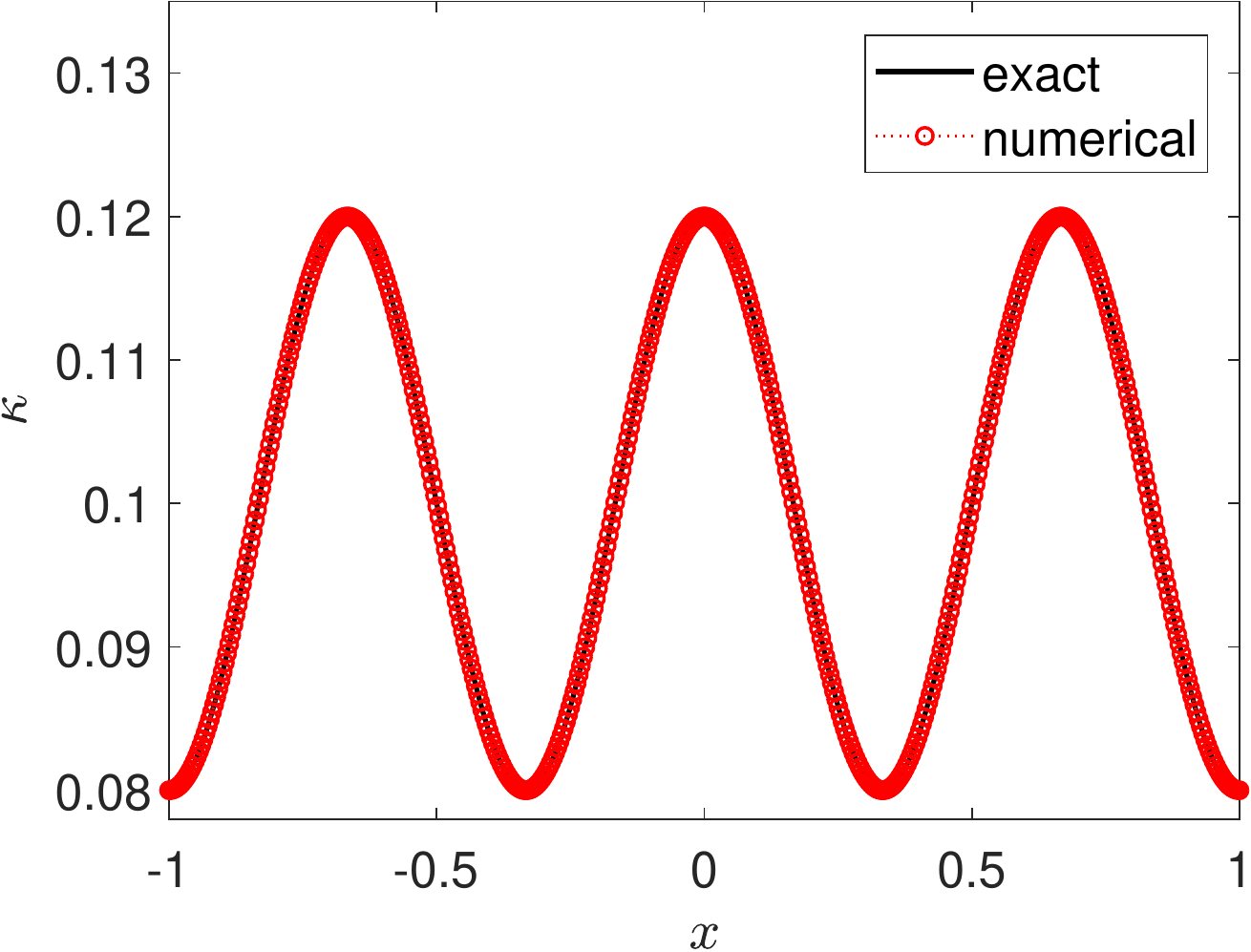}
		\caption{Example 4 with noiseless data and polynomial
                  order $n=40$. Left: recovery
                  of $\alpha$; Right: recovery of $\kappa $.}
	\label{ex4}
	\end{center}
\end{figure}

This nonlinear equation represents a departure from the linear
advection-diffusion equation discussed in the paper. 
Although our theoretical results do not apply here, 
the proposed numerical approaches still apply. We focus on noiseless
data case by Galerkin method. Similar to the other examples, data are
collected over 50 uniform grids in the temporal domain and 100 Gauss point
plus boundaries in the spatial domain.
Polynomial space of $\mathcal{P}_{1}^{40}$ is used as the
approximation and testing space for both $\alpha(x)$ and $\kappa(x)$.
The recovered results are shown in Fig.~\ref{ex4}. Good agreement
with the true $\alpha(x)$ and $\kappa(x)$ can be seen. The relative
errors are $9.7347\times 10^{-5}$ for $\alpha(x)$ and  $4.0632\times
10^{-5}$ for $\kappa(x)$.
%
%In this example, we use only a small time window $[0, 0.2]$ to collect data while avoiding potential shock. However, the recovery results turns out to be quite satisfactory. As we can see, when $n = 40$, we can reduce error to level $O(10^{-5})$. 
%
%
%%%%%%%%%%%%%%%%%%
%\begin{figure}[htbp]
%	\label{fig1_ex4}
%	\centering
%	\includegraphics[width=0.618\textwidth]{figures/ex4_burgers/ex4_burgers_solution_u_new-eps-converted-to.pdf}
%	\caption{Clean data $u$ at $t=0$ and $t=0.2$ for Example 4.}
%\end{figure}
%

\subsection{Example 5: 2D Advection-Diffusion Equation}

\begin{figure}[htbp]
	\begin{center}
	\includegraphics[width=0.48\textwidth]{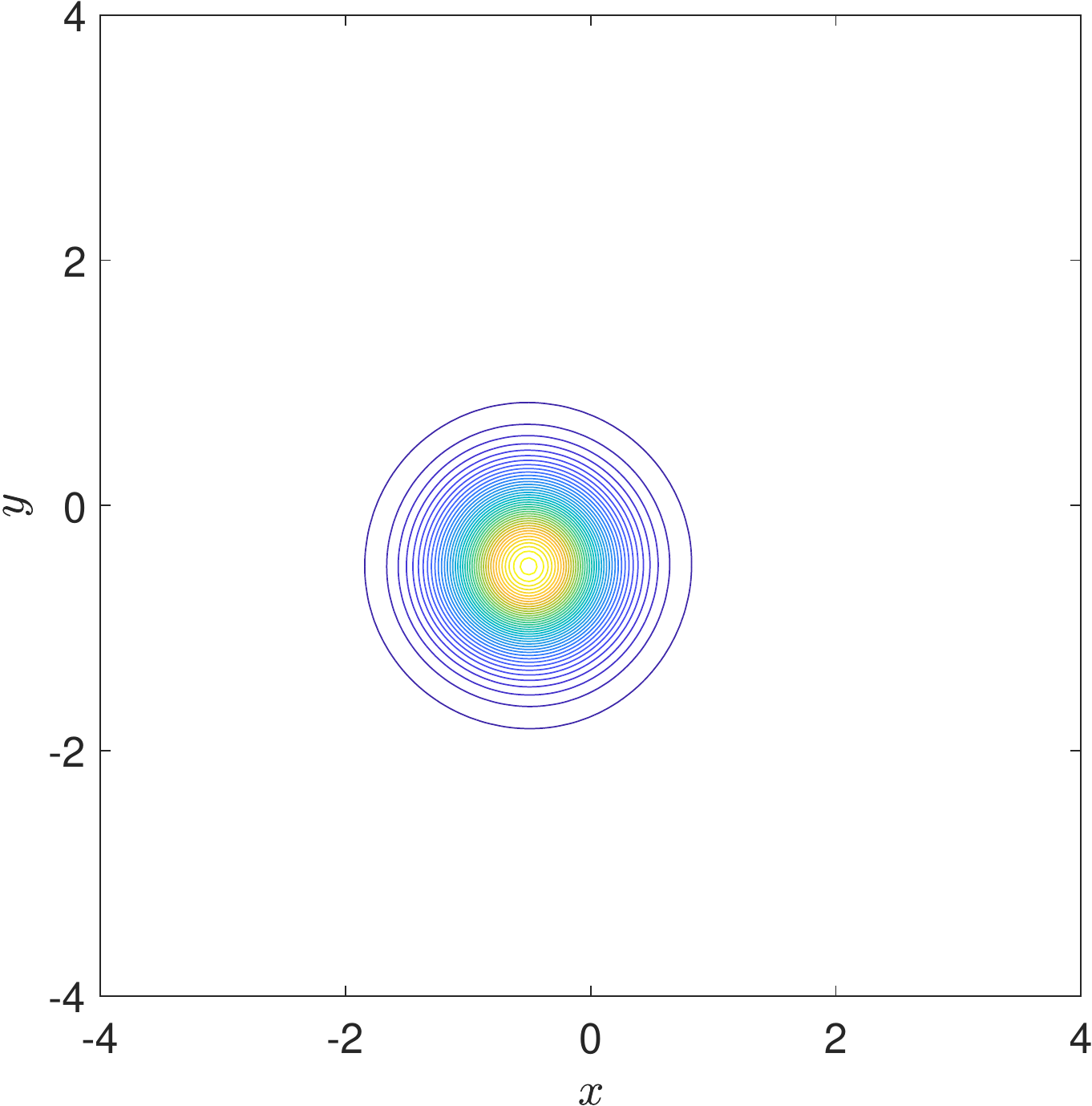}
	\includegraphics[width=0.48\textwidth]{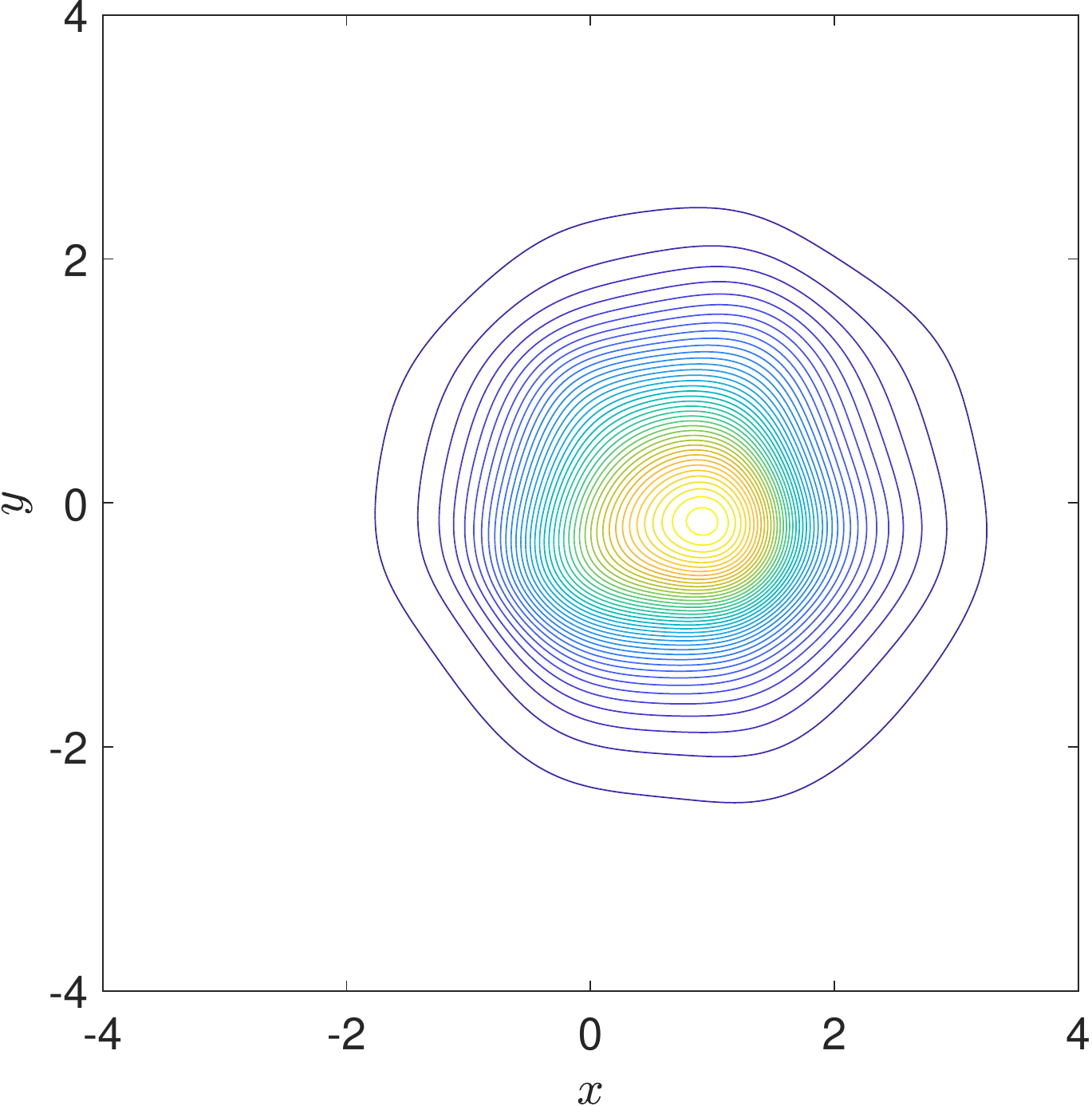}
	\caption{Example 5: State variable $u$ at Initial and final stage.}
	\label{ex5a}
\end{center}
\end{figure}

We finally consider a 2D advection-diffusion equation 
\begin{equation}
\frac{\partial u}{\partial t}(t,{\bm x}) = - \nabla \cdot ({\bm
  \alpha} u) + \nabla \cdot (\kappa \nabla u),
\qquad \textrm{in } (0,4] \times (-1,1)^2,
\end{equation}
where 
\begin{equation*}
\begin{split}
\bm{\alpha(x)} &= {\left (\widetilde{\alpha}_x y(1+\delta_{\alpha}
    \sin(\omega x)),-\widetilde{\alpha}_y x(1+\delta_{\alpha}
    \sin(\omega y)) \right)}^\top, \\
\kappa({\bm x})& = \widetilde{\kappa}\left(3 + \delta_{\kappa} \sin(\omega x) +
\delta_{\kappa} \cos(\omega y) \right),
\end{split}
\end{equation*}
with 
$\widetilde{\alpha}_x = 1$, $\widetilde{\alpha}_y = 1$,
$\delta_{\alpha} = 0.1$, $\widetilde{\kappa} = 0.02$, $\delta_{\kappa}
= 1$, and $\omega = \pi$. 
The initial condition is set as
$$
u(0,x,y) = \frac{1}{\sqrt{ (2\pi)^2 \sigma_x^2 \sigma_y^2}}{\rm
  e}^{\left(-\frac{1}{2}
  \frac{(x-\mu_x)^2}{\sigma_x^2}-\frac{1}{2}\frac{(y-\mu_y)^2}{\sigma_y^2}\right)},
$$
with
$\mu_x = \mu_y = -0.5$, $\sigma_x^2 =\sigma_y^2 = 0.2$.
The details of the numerical solver are in Table \ref{table_ex5}. The
solutions of the state variable at the initial and final time are
shown in Fig.~\ref{ex5a}, for demonstration purpose.
\begin{table}[htbp]
	\begin{center}
		\caption{PDE solver information for 2-D convection-diffusion equation in Example 5.}
		\begin{tabular}{ |c|c| }
			\hline
			data $u$ time domain  & [0,4]  \\  
			\hline
			data $u$ space domain  & ${[-4,4]}^2$   \\ 
			%\hline
			%initial condition & $u(0,x,y) = \frac{1}{\sqrt{ (2\pi)^2 \sigma_x^2 \sigma_y^2}}e^{(-\frac{1}{2} \frac{(x-\mu_x)^2}{\sigma_x^2}-\frac{1}{2}\frac{(y-\mu_y)^2}{\sigma_y^2})}, \mu_x = -0.5, \mu_y = -0.5, \sigma_x^2 = 0.2, \sigma_y^2 = 0.2 $\\
			\hline
			boundary condition & Dirichlet condition\\
			\hline 
			scheme in time  & Crank-Nicolson \\ 
			\hline
			time step size $\Delta t$ & $10^{-3}$ \\
			\hline
			scheme in space & Fourier collocation  \\
			\hline
			collocation number $N_{coll}$ & tensor points $80 \times 80$\\
			\hline
		\end{tabular}
	\label{table_ex5}
	\end{center}
\end{table}

To recover $\bm{\alpha(x)}$ and $\kappa({\bm x})$, we collection
solution data over  $200$ uniformly distributed grids in the temporal
domain and  $80 \times 80$ tensor Gauss points in the interior of the
spatial domain, along with $80$ Gauss points on each of the boundary
edges. We use $\mathcal{P}_2^8$ as approximation and testing space. 

\begin{figure}[htbp]
	\begin{center}
	\includegraphics[width=0.99\textwidth]{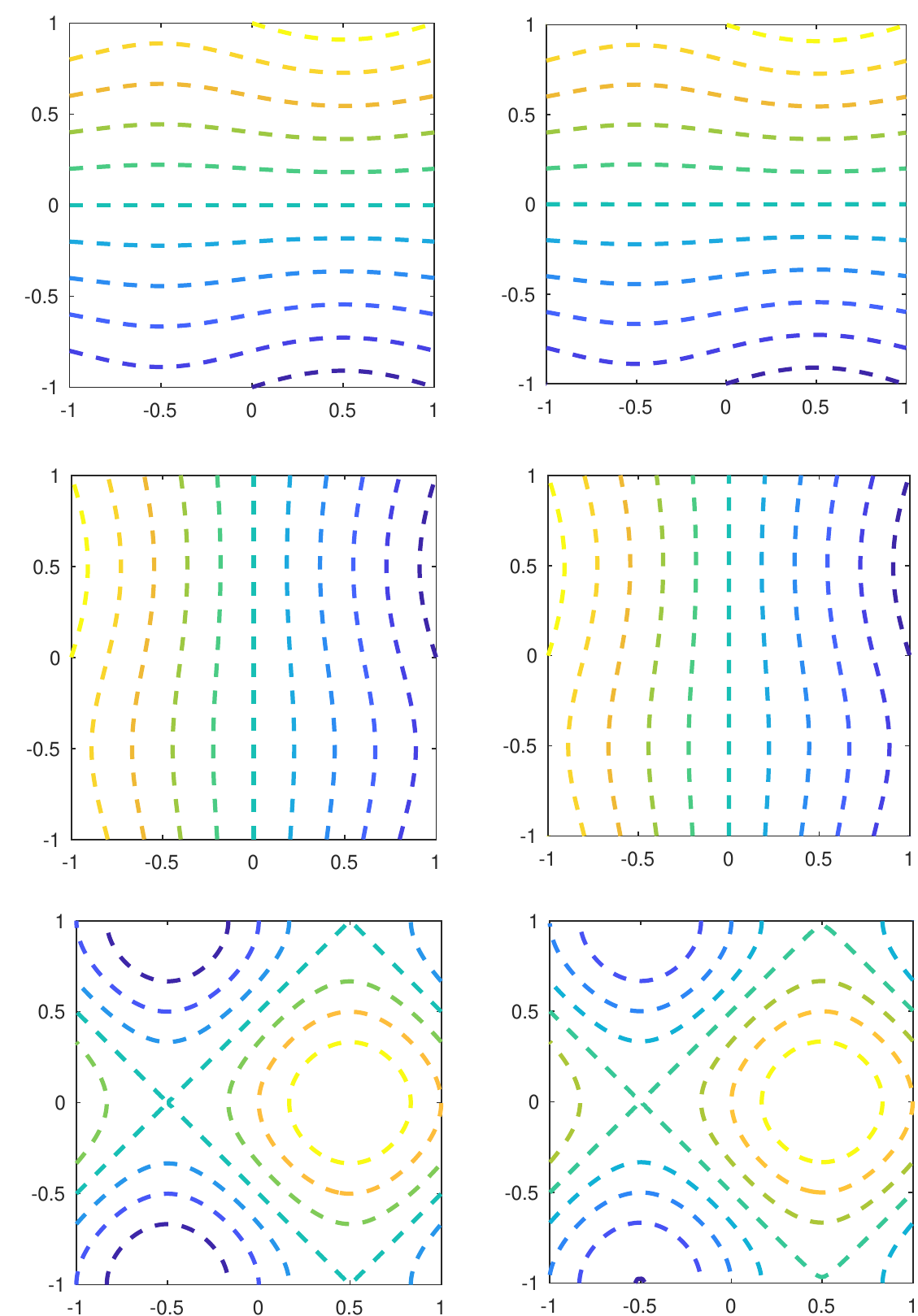}
	\caption{Example 5: Comparison of true (left column) and recovered
          (right column) parameter functions. From top to bottom
          $\alpha_1({\bm x})$,
          $\alpha_2({\bm x})$ and $\kappa({\bm x})$.}
	\label{fig5}
\end{center}
\end{figure}

The recovered results for $\bm{\alpha}(x) =
({\alpha}_1({\bm x}),{\alpha}_2({\bm x}))^{\top}$ and $\kappa({\bm x})$ are shown in
Fig.~\ref{fig5}, obtained via 
Galerkin method using noiseless data. Visual comparison with the true
functions shows good agreement. More detailed examination shows that
the relative errors in the recovered solutions are
$5.5235 \times 10^{-4} $ for $\alpha_1({\bm x})$,
$4.0274\times 10^{-4}$ for $\alpha_2({\bm x})$, and $6.9119\times 10^{-4}$
for $\kappa({\bm x})$. Results of noisy data case are not shown, as they are
visually similar to the noiseless case and with errors dominated by
the data noise.
%%%%%%%%%%%%%%%%%%

\section{Conclusion} \label{sec:conclusions}

In this paper, we studied the problem of identifying unknown parameter functions embedded in time-dependent partial differential equations (PDEs) using  observational data of the state variables. Using linear advection-diffusion type equations, we conducted theoretical analysis on the solvability of the problem and derived conditions under which unique recovery can be obtained. 
We then presented numerical approaches applicable for general PDEs. Two types of approaches, Galerkin and collocation, are presented.
While the collocation approach is straightforward to implement, the Galerkin method is preferred because its use of weak form avoids the use of much spatial derivatives of the state variables. In many practical cases when only data of the state variables are available, estimating derivatives
often induce additional errors, especially when data contain noises.

%\input Appendix
%\section*{Acknowledgment}
%This work is in part supported by AFOSR, DOE, NNSA, NSF.

\bibliographystyle{plain}
\bibliography{neural,LearningEqs,identification}

\end{document}